\documentclass[11pt]{article}
   \usepackage{amsfonts,amssymb}
\usepackage{latexsym,tikz}
\usepackage{amsmath}
\input xy
\xyoption{all}

\title{Multiple Dedekind Zeta Functions}
\author{Ivan Emilov Horozov}
\date{February 2, 2008}
\newcommand \nc {\newcommand}
\nc \proof {\noindent {\em{Proof.\/ }}} \nc \qed {$\Box$\hfill} 
\newtheorem{theorem}{Theorem}[section]
\newtheorem{lemma}[theorem]{Lemma}
\newtheorem{proposition}[theorem]{Proposition}
\newtheorem{corollary}[theorem]{Corollary}
\newtheorem{definition}[theorem]{Definition}
\newtheorem{example}[theorem]{Example}
\newtheorem{remark}[theorem]{Remark}
\newtheorem{conjecture}[theorem]{Conjecture}
\newtheorem{question}[theorem]{Question}
\nc \bth[1] {\begin{theorem}\label{t#1} } \nc \ble[1]
{\begin{lemma}\label{l#1} } \nc \bpr[1]
{\begin{proposition}\label{p#1} } \nc \bco[1]
{\begin{corollary}\label{c#1} } \nc \bde[1]
{\begin{definition}\label{d#1}\rm } \nc \bex[1]
{\begin{example}\label{e#1}\rm } \nc \bre[1]
{\begin{remark}\label{r#1}\rm } \nc \bcon[1]
{\begin{conjecture}\label{con#1}\rm } \nc \bque[1]
{\begin{question}\label{que#1}\rm }
\def \P {{\mathbb P}}

\def \N {{\mathbb N}}
\def \Z {{\mathbb Z}}
\def \Q {{\mathbb Q}}
\def \R {{\mathbb R}}
\def \C {{\mathbb C}}
\def \H {{\mathbb H}}

\def \Im { {\mathrm{Im}} }

\begin{document}

%%%%%%%%%%%%%%%%%%%%%%%%%%%%%%%%%%%%%%%%%%%%%%%%%%%%%%%%%%%%%%%%%%%%%%%
%%%%%%
\title{{\LARGE\bf{Non-commutative Hilbert Modular Symbols}}}

\author{
Ivan ~Horozov
\thanks{E-mail: horozov@math.wustl.edu}
\\ \hfill\\ \normalsize \textit{Department of Mathematics,}\\
\normalsize \textit{Washington University in St. Louis,}\\
\normalsize \textit{One Brookings Drive, Campus Box 1146}\\
\normalsize \textit {Saint Louis, MO 63130, USA }  \\
}
\date{}
\maketitle

\begin{abstract}
The main goal of this paper is to construct non-commutative Hilbert modular symbols. However, we also construct commutative Hilbert modular symbols. Both the commutative and the non-commutative Hilbert modular symbols are generalizations of Manin's classical and non-commutative modular symbols.
We prove that many cases of (non-)commutative Hilbert modular symbols are periods in the sense on Kontsevich-Zagier. Hecke operators act naturally on them. 

Manin defines the non-commutative modilar symbol in terms of iterated path integrals.
In order to define non-commutative Hilbert modular symbols, we use a generalization of iterated path integrals to higher dimensions, which we call iterated integrals on membranes. Manin examines similarities between non-commutative modular symbol and multiple zeta values both in terms of infinite series and in terms of iterated path integrals. Here we examine similarities in the formulas for non-commutative Hilbert modular symbol and multiple Dedekind zeta values, recently defined by the author, \cite{MDZF}, both in terms of infinite series and in terms of iterated integrals on membranes.
\end{abstract}
%%%%%%%%%%%%%%%%%%%%%%%%%%%%%%%%%%%


\tableofcontents
\setcounter{section}{-1}
%%%%%%%%%%%

\section{Introduction}
\label{Intro}
%    \begin{multline}
%          \sin A \cos B = \frac{1}{2}\left[ \sin(A-B)+\sin(A+B) \right] \\
%          \sin A \sin B = \frac{1}{2}\left[ \sin(A-B)-\cos(A+B) \right] \\
%          \cos A \cos B = \frac{1}{2}\left[ \cos(A-B)+\cos(A+B) \right] \\
%          \end{multline}

%\begin{align}
%&\zeta_{K;C}(s_1,\dots,s_1;\dots;s_m,\dots,s_m)=\\
%&\sum_{\alpha_1,\dots,\alpha_m\in C}
%\frac{1}{N(\alpha_1)^{s_1}N(\alpha_1+\alpha_2)^{s_2}\cdots N(\alpha_1+\dots+\alpha_m)^{s_m}}.
%\end{align}

Classical elliptic modular symbols were introduced by Birch \cite{Bi} and Manin \cite{Manin} in connection with Birch - Swinnerton-Dyer conjecture for certain congruence subgroups of $SL_2(\Z)$. We recall that a modular symbol $\{p,q\}$ is 
associated to a pair of cusp points $p,q\in \P^1(\Q)$ 
on the completed upper half plane $\H^1\cup\P^1(\Q)$. 
One can think of the modular symbol $\{p,q\}$ as a homology class of the geodesic connecting $p$ 
and $q$, in $H_1(X_\Gamma,\{cusps\})$, where $X_\Gamma$ is the modular curve associated to a 
congruence subgroup of $SL_2(\Z)$. One can pair $\{p,q\}$ with a cusp form $f$ by
\[
\{p,q\}\times f\mapsto\int_p^qfdz
\]
If $f$ is a cusp form of weight $2$ 
then $fdz$ can be considered as cohomology class in $H^1(X_\Gamma)$.
This gives a pairing between homology (Betti) and cohomology (de Rham) that leads to periods. Modular symbols are a useful tool applied to $L$-functions and computation of cohomology groups. For a review of such topics, one can consult \cite{M}.

Their theory was developed by Manin, Drinfeld,  Shokurov, Mazur, \cite{Manin, Dr, Shok, Maz}.    Later the theory was extended to higher ranks by A. Ash, Rudolf, A. Borel, Gunnels \cite{AR, AB, Gu}.

Elliptic modular symbols are important tool in the study of modular forms. They are particularly useful in computations with modular forms. J. Cremona  designed  algorithms for computations with elliptic curves, based on modular symbols ("modular symbol algorithm"), see \cite{Cre}.   Some of the applications include computations of homology and cohomology.  Also, the study of special  values of $L$-functions became a vast area of applications of classical modular symbols, see \cite{MS, KMS}.

Later   W. Stein also has contributed to the difficult area of computations with modular forms. See his excellent book  \cite{Ste}, which contains both theory and computational methods. For higher rank groups one can consult the Appendix by P. Gunnels in the same book.

Manin's non-commutative modular symbol \cite{Man} is a generalization of both the classical modular symbol and of multiple zeta values in terms of Chen's iterated integral theory in the holomorphic setting. Manin shows that the non-commutative modular symbol is a non-commutative $1$-coclycle.
He  also shows that each of the iterated integrals on Hecke eigenforms that enter in the non-commutative modular symbol are periods.  

The main goal of this paper is to construct non-commutative Hilbert modular symbols. However, we also construct an analog of the classical modular symbol for Hilbert modular varieties. Both symbols are generalizations of the corresponding Manin's constructions. 

We compute explicit integrals in terms of the non-commutative Hilbert modular symbol of type {\bf{b}} and present similar formulas for the recently defined multiple Dedekind zeta values (see \cite{MDZF}). We prove that the iterated integrals on membranes that enter in the non-commutative modular symbol of type {\bf{c}} are periods. We also give some explicit and some categorical arguments in support of a conjecture that a certain type of non-commutative Hilbert modular symbol satisfies a non-commutative $2$-cocycle condition.

Before describing our results let us recall the non-commutative modular symbol of Manin \cite{Man}. 
Let $\nabla=d-\sum_{i=1}^mX_if_idz$ 
be a connection on the upper half plane, where $f_1,\dots,f_m$ 
are cusp forms and $X_1,\dots,X_m$ are formal variables. 
One can think of $X_1,\dots,X_m$ as constant square matrices of the same size.

Let $J^b_a$ be the parallel transport of the identity matrix $1$ at the point $a$ to the matrix $J^b_a$ 
at the point $b$. Alternatively, $J^a_b$ can be written as a generating series of iterated path integrals
of the forms $f_1dz,\dots,f_mdz$, (see \cite{Ch} and \cite{Man}), namely,
\[
J^b_a
=
1
+
\sum_{i=1}^mX_i\int_a^bf_idz
+
\sum_{i,j=1}^m
X_iX_j\int_a^b f_idz\cdot f_jdz
+\cdots
\]
Then $J_a^bJ_b^c=J_a^c$. 
This property leads to the $1$-cocycle $c_a^1(\gamma)=J^a_{\gamma a},$
which is the non-commutative modular symbol (see \cite{Man} and Section 1 of this paper).
If $f_1,\dots,f_m$ are normalized cusp Hecke eigenforms then each iterated integral appearing in the generating series $J_a^b$ is a period.
In this paper we introduce   both commutative and non-commutative modular symbols for Hilbert modular surfaces. As it turned out we need some  new tools in comparison to the classical modular symbols.  In particular for the non-commutative Hilbert modular symbol we need iterated integrals in dimension higher than one. We introduce them and study their properties in the special case of Hilbert modular surfaces.
In this paper, we construct both commutative and non-commutative modular symbols for the Hilbert modular group $SL_2({\cal{O}}_K)$. For the Hilbert modular group, 
one may consult \cite{B} and \cite{F}. In the case of Hilbert modular surface, it is not possible to repeat Manin's constructions for the  non-commutative modular symbols, since  the integration domain is two-dimensional over the complex numbers. Instead, we develop a new approach (Section 2), which we call iterated integrals on membranes. This is a higher dimensional analogue of iterated path integrals.
In Subsection 3.7, we explore similar relations between non-commutative Hilbert modular symbols and 
multiple Dedekind zeta values  (see \cite{MDZF}).

In Section 3, we associate modular symbols for $SL_2({\cal{O}}_K)$ to geodesic triangles and 
geodesic diangles ($2$-cells whose boundary has two vertices and two edges, which are geodesics.) 
We are going to explain how the geodesic triangles and the geodesic diangles are constructed.
Consider $4$ cusp points in $\H^2\cup \P^1(K)$.
We can map every three of them to $0$, $1$ and $\infty$ with a linear fractional transformation 
$\gamma\in GL_2(K)$. There is a diagonal map $\H^1\rightarrow \H^2$, whose image $\Delta$ 
contains $0$, $1$ and $\infty$. We can take a pull-back of $\Delta$ with respect to the map $\gamma$ in order to obtain a holomorphic  (or anti-holomorphic) curve that passes through the given three points. If $\det\gamma$ is totally positive or totally negative then $\gamma^*\Delta$ 
is a holomorphic curve in $\H^2$.
If $\det \gamma$ is not totally positive or totally negative (that is in one of the real embedding it is 
positive and in the other it is negative) then $\gamma^*\Delta$ is anti-holomorphic. This means that it 
is holomorphic in $\H^2$ if we conjugate the complex structure in one of the copies of $\H^1$. The 
same type of change of the complex structure is considered in \cite{F}.

On each holomorphic, or anti-holomorphic curve $\gamma^*\Delta$, there is a unique geodesic 
triangle connecting the three given points. However, if we take two of the points, we see that they 
belong to two geodesic triangles. Thus they belong to two holomorphic, (or anti-holomorphic) curves. 
Therefore, there are two geodesic connecting the two points - each lying on different holomorphic  
(or anti-holomorphic) curves, as faces of the corresponding geodesic triangles, defining the curves. 
There are two pairings that we consider: the first one is an integral of a cusp form over a geodesic 
triangle and the second one is an integral of a cusp form over a geodesic diangle. If we integrate a 
holomorphic 2-form coming from a  cusp form over a 
geodesic triangle, we obtain $0$, if the triangle lies on an holomorphic curve. Thus the only non-zero 
pairings come from integration of a cusp form over a diangle or over a triangle, lying on an 
anti-holomorphic curve.

Now let us look again at the four cusp points together with the geodesics that we have just  described. We obtain four geodesic triangles, corresponding to each triple of points among the four points, and six diangles, corresponding to the six ``edges" of a tetrahedron with vertices the four given points. 
Thus, we obtain a ``tetrahedron" with thickened edges. We will use tetrahedrons with thickened edges as an intuition for a non-commutative $2$-coclycle relation (see Conjecture \ref{thm c2}) for the non-commutative Hilbert modular symbol, which is an analogue of Manin's non-commutative 1-cocycle relation for the non-commutative modular symbol.

Usually, the four vertices are treated as a tetrahedron and a $2$-cocycle is functional on the faces, considered as 2-chains. 
The boundary is defined as a sum of the $2$-cocycles on each of the faces (which are triangles). The 
boundary of the tetrahedron gives a boundary relation of a $2$-cocycle.

In our case the analogue of a $2$-cocycle is a functional on diangles and on triangles. And the 
boundary map is a sum over the faces of the thickened tetrahedron. Thus, the faces of the thickened 
tetrahedron are four triangles and six diangles, corresponding to the six edges of a tetrahedron.

%In this paper we construct a non-commutative Hilbert modular symbol, which is a generalization to 
%Manin's non-commutative modular symbol for congruence subgroups of $SL_2(\Z)$. Manin used 
%iterated path integrals in the upper half plane. In contrast, we define here a suitable higher 
%dimensional analogue of iterated path integrals, which we call iterated integrals over membranes.
%Manin proved that the non-commutative modular symbol is a non-abelian first cocycle. 
%Here, we prove analogues statement. We show that the non-abelian Hilbert modular symbol 
%satisfies a relation resembling a non-commutative 2-cocycle relation. In particular it gives a relation 
%among commutative Hiblert modular symbols that we define here.

We show that the geodesics on the boundary of a diangle or of a geodesic triangle lie on a 
holomorphic curves $\gamma^*\Delta$ for various elements $\gamma$ with totally positive or totally 
negative determinant. This implies that when we take the quotient by a Hilbert modular group the 
holomorphic curve $\gamma^*\Delta$ becomes Hirzebruch-Zagier divisor \cite{HZ}. Then we prove 
that the commutative Hilbert modular symbols paired with a cusp forms of weight $(2,2)$ give periods 
in the sense of \cite{KZ}.

In order to construct a non-commutative Hilbert modular symbol, first we define a suitable 
generalization of iterated path integrals,
which we call iterated integrals on membranes (see Section 2). 
We choose the word ``membrane" since such integrals are invariant under suitable variation of the 
domain of integration. 

There is a topological reason for considering non-commutative Hilbert modular symbol as opposed to 
only commutative one. 
Let us first make such comparison for the case of $SL_2(\Z)$. The commutative modular symbol 
captures $H_1(X_\Gamma)$, while the non-commutative symbol captures the rational homotopy type 
of the modular curve $X_\Gamma$. Now, let $\tilde{X}$ be a smooth Hilbert modular surface, 
by which we 
mean the minimal desingularization of the Borel-Baily compactification due to Hirzebruch. 
Then the rational fundamental group of a Hilbert modular surface vanishes, $\pi_1(\tilde{X})_\Q=0$, 
(see \cite{B}). The non-commutative Hilbert modular symbol is an attempt to capture more from the 
rational homotopy type compared to what $H_2(\tilde{X})$ captures.

For convenience of the reader, first we define type {\bf{a}} iterated integrals on membranes 
(Definition \ref{def membranes 2a}). They are 
simpler to define and more intuitive. However, they do not have enough properties. 
(For example, they do not have an integral shuffle relation.) 
Then we define the type {\bf{b}} iterated integrals on membranes 
(Definition \ref{def membranes 1b}), which involves 
two permutations. Type {\bf{b}} has integral shuffle relation (Theorem \ref{thm shuffle} (i)) 
and type {\bf{a}} is a particular case of type {\bf{b}}.

We are mostly interested in iterated integrals of type {\bf{b}}. If there is no index specifying the type of 
iterated integral over membranes, we assume that it is of type {\bf{b}}.

Similarly to Manin's approach, we define a generating series of iterated integrals over membrane of 
type {\bf{b}} over $U$, which we denote by $J(U)$.
We also define a shuffle product of generating series of iterated integrals over membranes of type 
{\bf{b}} (see Theorem \ref{thm shuffle} 
part (iii)),
\[\phi(J(U_1)\times_{Sh}J(U_2))=J(U_1\cup U_2)\] 
for disjoint manifolds with corners of dimension $2$,
$U_1$ and $U_2$, as subsets of $\H^2\cup\P^1(K)$, (see \cite{BS}). 
This shuffle product generalizes the composition of generating series of iterated path integrals, namely, 
$J_a^bJ_b^c=J_a^c$, to dimension $2$. Note that similar definition is also possible in higher 
dimensions. Also $J(U)$ is invariant under homotopy. This allows to 
consider cocycles and coboundaries, where the relations use homotopy invariance and the values at 
the different cells can be composed via the shuffle product.

We define non-commutative Hilbert modular symbols which we call $c^1$ and $c^2$. Then $c^1$ is 
the functional $J$ on certain geodesic diangles and $c^2$ is the functional $J$ on geodesic triangles. 
Conjecturally, $c^1$ is a $1$-cocycle such that 
if we change the base point of $c^1$ then $c^1$ is modified by a coboundary. 
Also conjecturally, $c^2$ is a $2$-cocycle up to finitely many multiples of different values of $c^1$. 
Also, if we change the base point of $c^2$ then $c^2$ is modified by a coboundary up to a finitely 
many multiples of different values of $c^1$. 
In Subsection 3.5 we give explicit formulas in support of the interpretation of the non-commutative symbols as co-cycles. 

In Subsection 3.6, we give a categorical construction, which might help proving that the non-commutative symbols as co-cycles.

In Subsection 3.7, we define multiple $L$-values associated to cusp forms and we compare 
them to multiple Dedekind zeta values (see \cite{MDZF}).

We also briefly recall the Riemann zeta values and multiple zeta values (MZVs). The Riemann zeta values are defined as
$$\zeta(k)=\sum_{n>0}\frac{1}{n^k},$$
where $n$ is an integer.
MZVs are defined as
$$\zeta(k_1,\dots,k_m)=\sum_{0<n_1<\dots<n_m}\frac{1}{n_1^{k_1}\dots n_m^{k_m}},$$
where $n_1,\dots,n_m$ are integers. The above MZV is of depth $m$.
Riemann zeta values $\zeta(k)$ and MZV
$\zeta(k_1,\dots,k_m)$ were defined by Euler \cite{Eu} for $m=1,2$. 

The common feature of MZVs and the non-commutative modular symbol is that they both can be 
written as iterated path integrals (see \cite{G1}, \cite{G2}). Moreover, Manin's non-commutative modular symbol resembles the generating series of MZV, which is 
the Drinfeld associator. Let us recall that the Drinfeld associator is a generating series of iterated integrals of the type $J_a^b$ associated to the connection
\[\nabla=d-A\frac{dx}{x}-B\frac{dx}{1-x}\]
on $Y_{\Gamma(2)}=\P^1-\{0,1,\infty\}$. One can think of $Y_{\Gamma(2)}$ as the modular curve 
associated to the congruence subgroup $\Gamma(2)$ of $SL_2(\Z)$. 
Then the differential forms $\frac{dx}{x}$ and $\frac{dx}{1-x}$ 
are Eisenstein series of weight $2$ on the modular curve $Y_{\Gamma(2)}$.

Relations between MZV and modular forms have been examined by many authors.
For example, Goncharov has considered a mysterious relation between MZV (multiple zeta values) of 
given weight and depth $3$ and cohomology of $GL_3(\Z)$ (see \cite{G2} and \cite{G3}), which is 
closely related to the cohomology of $SL_3(\Z)$. In pursue for such a relation in depth 4, Goncharov 
has suggested and the author has computed the group cohomology of $GL_4(\Z)$ with coefficients in 
a family of representations, \cite{GL4}.  Another relation between modular forms and MZV is presented 
in \cite{ZGK}.

Similarly to Manin's approach, we explore relations between the non-commutative Hilbert modular symbols and multiple Dedekind zeta values (see \cite{MDZF}). Let us recall multiple Dedekind zeta values. Let each of $C_1,\dots,C_m$ be a suitable subset of the ring of integers ${\cal{O}}_K$ of a number field $K$. We call each of $C_1,\dots,C_m$ a cone. Then
{\bf{multiple Dedekind zeta values}} are defined as
\begin{equation}
\nonumber
\zeta_{K;C_1,\dots,C_m}(k_1,\dots,k_m)
=
\sum_{\alpha_i\in C_i\text{ for }i=1,\dots,m}
\frac{1}{N(\alpha_1)^{k_1}N(\alpha_1+\alpha_2)^{k_2}\cdots N(\alpha_1+\cdots+\alpha_m)^{k_m}}.
\end{equation}
The connection between non-commutative Hilbert modular symbols and multiple Dedekind zeta  values is both in similarities in the infinite sum formulas and in the definition in terms of iterated integrals on membranes (see \cite{MDZF}).

We consider a non-commutative Hilbert modular symbol of type {\bf{b}} over one diangle and 
compare it with (multiple) Dedekind zeta values with summation over one discrete cone \cite{MDZF}. 
However, in this case the two series look very different. We obtain that the multiple $L$-values are non-
commutative modular symbols defined as $J$ evaluated at an infinite union of diangles. We obtain that 
such $L$-values are very similar to the sum of multiple Dedekind zeta values, in the same way as the integrals in the Manin's non-commutative modular symbol are similar to the multiple zeta values (MZV). Then the sum of the multiple Dedekind zeta values is over an infinite union of cones. 
The idea to consider cones originated by Zagier \cite{Zagier} and more generally by Shintani \cite{Sh}.

Classical or commutative modular symbols for $SL_3(\Z)$ and $SL_4(\Z)$ were constructed by Ash and Borel \cite{AB} and Ash 
and Gunnels \cite{AG}. For $GL_2({\cal{O}}_K)$, where $K$ is a real quadratic field, Gunnels and 
Yasaki  have defined a modular symbol
based on Voronoi decomposition of a fundamental domain, in order to compute the $3$-rd cohomology 
group of $GL_2({\cal{O}}_K)$, (see \cite{GY}). (For the Hilbert modular group, $SL_2({\cal{O}}_K)$, 
one may consult \cite{B} and \cite{F}.) In contrast, here we use a geodesic triangulation of 
$\H^2/SL_2({\cal{O}}_K)$. We are interested mostly in  $2$-cells, whose boundaries are geodesics.
One of the (commutative) symbols that we define here, resembles combinatorially the symplectic 
modular symbol of Gunnells, \cite{Gu}. However, the meanings of the two types of symbols and their 
approaches are different.

There are several different directions for further work on Hilbert modular symbols. First of all, the commutative Hilbert modular symbols have good behavior with respect to Hecke operators. It will be interesting to extend the Hecke operators to the cases of higher equal weight case $(k,k)$. To apply Hecke operators to the Hilbert modular groups one either assumes a trivial narrow class group or one has to work with adeles. Another possible continuation of the coronet work is to extend commutative Hilbert modular symbols to the adelic setting. Then, one may try to extend these properties  - higher equal weight cusp forms and Hecke operators in the adelic setting to the non-commutative Hilbert modular symbols.
Hopefully, the abelian Hilbert modular symbol would lead to computational tools for cohomology of some Hilbert modular groups with coefficients in various representations. 

For the non-commutative Hilbert modular symbols we expect that some of the continuations would be establishing the $2$-categorical framework that define non-abelian $2$-cohomology set. Such a task would also have applications to non-commutative reciprocity laws on algebraic surfaces. In dimension $1$, we have a non-commutative reciprocity law as a reciprocity law for a generating series of iterated path integrals on a complex curve \cite{rec1}. In dimension $2$ we have proven both the Parshin reciprocity and a new reciprocity for a $4$-function local symbols \cite{rec2} defined by the author, which are particular cases in the generating series. A $2$-categorical $2$-nd cohomology set would capture algebraically the generating series  of iterated integrals on membranes needed for the general reciprocity on algebraic surfaces.

Finally, we expect that the (non-)commutative Hilbert modular symbols would be used for studying $L$-functions and multiple $L$-functions together with their special values.

\section{Manin's non-commutative modular symbol}
In this section we would like to recall the definition and the main properties of the Manin's 
non-commutative modular symbol, (see \cite{Man}). In his paper \cite{Man}, Manin uses iterated path integrals on a modular curve and on its universal cover - the upper half plane. Our main constructions are parallel to some extend to Manin's approach and for that reason we recall it below. However, instead of iterated path integrals we introduce a new tool - iterated integrals on membranes (see Section 2). Only this notion is adequate for studying non-commutative Hilbert modular symbols by generalizing the iteration process to higher dimensions.

\subsection{Iterated path integrals}
Here we recall iterated path integrals (see also \cite{Parshin}, \cite{Ch}, \cite{G1}, \cite{Man}). In the next Section, we are going to generalize them to iterated integrals over membranes.

\begin{definition}
\label{def it path int}
Let $\omega_1,\dots,\omega_m$ be $m$ holomorphic $1$-forms 
on the upper half plane together with the cusps, $\H^1\cup \P^1(\Q)$. Let
\[g:[0,1]\rightarrow \H^1\cup \P^1(\Q),\]
be a piece-wise smooth path.
We define an iterated integral
\[\int_g\omega_1\dots\omega_m
=
\int\dots\int_{0<t_1<t_2\dots<t_m<1}g^*\omega_1(t_1)\wedge\dots\wedge g^*\omega_m(t_m).\]
\end{definition}
Let $X_1,\dots,X_m$ be formal variables. Consider the differential equation
\begin{equation}
\label{eq diff eq}
dF(\Omega)=F(\Omega)(X_1\omega_1+\dots X_n\omega_m)
\end{equation}
with values in the 
associative but non-commutative
ring of formal power series in the non-commuting variables $X_1,\dots,X_m$ over the ring of 
holomorphic functions on the upper half plane. There is a unique solution with initial condition for $F(\Omega)(g(0))=1$, at the starting point $g(0)$, equal to $1$.
Then $F(\Omega)$ at the end of the path, that is at the point $g(1)$, has the value
\begin{equation}
\label{eq solution}
F_g(\Omega)
=
1
+\sum_{i=1}^mX_i\int_g\omega_i
+\sum_{i,j=1}^mX_iX_j\int_g\omega_i\omega_j
+\sum_{i,j,k=1}^mX_iX_jX_k\int_g\omega_i\omega_j\omega_k+\dots
\end{equation}
Using the Solution \eqref{eq solution}  to Equation \eqref{eq diff eq}, we prove the following theorem.
\begin{theorem}
\label{thm composition}
Let $g_1$ and $g_2$ be two paths such that the end of $g_1$, $g_1(1)$ is equal 
to the beginning of $g_2$, $g_2(0)$. Let $g_1g_2$ denote the concatenation of $g_1$ and $g_2$. 
Then
\[F_{g_1g_2}(\Omega)=F_{g_1}(\Omega)F_{g_2}(\Omega).\]
\end{theorem}
\proof 
The left hand side is the value of the solution of 
the linear first order ordinary differential equation at the point $g_2(1)$. 
From the uniqueness of the solution, we have that the solution along $g_2$ gives the same result, 
when the initial condition at $g_2(0)$ is $F_{g_1}(\Omega)$. 
That result is $F_{g_1}(\Omega)F_{g_2}(\Omega)$.

The same result can be proven via product formula for iterated integrals. 
We need this alternative proof in order to use it for generalization to higher dimensions.
\begin{lemma} 
\label{lemma product of paths}(Product Formula) Let $\omega_1,\dots,\omega_m$ be holomorphic $1$-forms on $\C$ and $g_1,g_2$ be two paths such that the end of $g_1$ is the beginning of $g_2$, that is $g_1(1)=g_2(0)$. As before we denote by $g_1g_2$ the concatenation of the paths $g_1$ and $g_2$. Then
\[\int_{g_1g_2}\omega_1\cdots\omega_m
=
\sum_{i=0}^m\int_{g_1}\omega_1\cdots\omega_i
\int_{g_2}\omega_{i+1}\cdots\omega_m.\]
\end{lemma}
\proof Let $g_1:[0,1]\rightarrow \C$ and let $g_2:[1,2]\rightarrow \C$. We consider the concatenation $g_1g_2$ as a map $g_1g_2:[0,2]\rightarrow \C$ such that its restriction to the interval $[0,1]$ gives the path $g_1$ and its restriction to the interval $[1,2]$ gives $g_2$.
From Definition \ref{def it path int}, we have that
\[\int_{g_1g_2}\omega_1\cdots\omega_m=\int\cdots\int_{0<t_1<\cdots<t_m<2}(g_1g_2)^*\omega_1(t_1)\wedge\cdots\wedge (g_1g_2)^*\omega_m(t_m).\]
In the domain of integration 
$0<t_1<\cdots<t_m<2$ insert the number $1$. Geometrically, we cut the simplex $0<t_1<\cdots<t_m<2$ into a disjoint union of products of pairs of simplexes such that its $t_k\in[0,1]$ for $k\leq i$ and $t_k\in [1,2]$ for $k>i$. Thus, the union is over distinct values of $i$ for $i=0,\dots,m$. And for each fixed $i$ the two simplexes are $0<t_1<\cdots<t_i<1$ and $1<t_{i+1}<\cdots<t_m<2.$
Then we have 
\begin{align*}
&\int_{g_1g_2}\omega_1\cdots\omega_m=\\
=&\sum_{i=0}^n\int\cdots\int_{0<t_1<\cdots<t_i<1;\,1<t_{i+1}<\cdots<t_m<2}
 (g_1g_2)^*\omega_1(t_1)\wedge\cdots\wedge (g_1g_2)^*\omega_m(t_m)
=\\
=&\sum_{i=0}^n\left(\int\cdots\int_{0<t_1<\cdots<t_i<1}
 g_1^*\omega_1(t_1)\wedge\cdots\wedge g_1^*\omega_i(t_i)
 \right)\times\\
 &\times\left(
\int\cdots\int_{1<t_{i+1}<\cdots<t_m<2}
 g_2^*\omega_{i+1}(t_{i+1})\wedge\cdots\wedge g_2^*\omega_m(t_m)\right)=\\
 =&\sum_{i=0}^m\int_{g_1}\omega_1\cdots\omega_i
\int_{g_2}\omega_{i+1}\cdots\omega_m.
\end{align*}

\begin{definition}
\label{def shuffle of sets} The set of all shuffles
$sh(i,j)$  is  a subset of all permutations of the set $\{1,2,\dots,i+j\}$ such that
\[\rho(1)<\cdots<\rho(i)\] 
and
\[\rho(i+1)<\cdots<\rho(i+j).\]
Such a permutation $\sigma$  is called a shuffle.
\end{definition}

\begin{lemma}
\label{lemma shuffle}
 (Shuffle Relation) Let $\omega_1,\dots,\omega_m$ be holomorphic $1$-forms on $\C$ and let $g$ be a path. Then
\[\int_g\omega_1\cdots\omega_i
\int_{g}\omega_{i+1}\cdots\omega_m=
\sum_{\sigma\in sh(i,m-i)}\int_g\omega_{\rho(1)}\cdots\omega_{\rho(m)},\]
where $sh(i,j)$ is the set of shuffles from Definition \ref{def shuffle of sets}.
\end{lemma}

\subsection{Manin's non-commutative modular symbol}
Now let $g$ be a geodesic connecting two cusps $a$ and $b$ 
in the completed upper half plane $\H^1\cup \P^1(\Q)$. Let 
$\Omega=\{f_1dz,\dots,f_mdz\}$ be a finite set of holomorphic forms with respect to a congruence subgroup $\Gamma$ of $SL_2(\Z)$, such that $f_1,\dots,f_m$ are cusp forms of weight $2$.
Let
\[J_a^b=F_g(\Omega).\]
As a reformulation of Theorem \ref{thm composition}, we obtain the following.
 \begin{lemma}
 \[J_a^bJ_b^c=J_a^c.\]
 \end{lemma}
 We give a direct consequence of it. 
 \begin{corollary}
\[J_b^a=(J_a^b)^{-1}.\]
\end{corollary}

Now we are ready to define Manin's non-commutative modular symbol.
Note that there is a natural action of $\Gamma$ on $J_a^b$. If $\gamma\in\Gamma$ then
$\gamma J_a^b$ is defined as $J_{\gamma a}^{\gamma b}$. If $f_1,\dots,f_m$ are cusp forms of weight $2$. Then $\omega_1=f_1dz,\dots,\omega_m=f_mdz$ are forms of weight $0$, that is, they are invariant forms with respect to the group $\Gamma$. Then 
\[\gamma J_a^b=F_{\gamma g}(\omega_1,\dots,\omega_m)=F_g(g^*\omega_1,\dots,g^*\omega_m)
=F_g(\omega_1,\dots,\omega_m)
=J_a^b.\]
Let $\Pi$ be a subgroup of the invertible elements $\C\ll X_1,\dots,X_m\gg$ with constant term $1$.
We extend action of $\Gamma$ on $J_a^b$ to a trivial action of $\Gamma$ on $\Pi$.

Following Manin, we present the key Theorem for and Definition of the non-commutative modular symbol.
\begin{theorem} Put
\[c^1_a(\gamma)=J^a_{\gamma a}.\]
Then $c^1_a$ represent a cohomology class in $H^1(\Gamma,\Pi)$ independent of the base point $a$
\end{theorem}
\proof First $c^1_a$ is a cocycle:
\[dc^1_a(\beta,\gamma)=J^a_{\beta a}(\beta\cdot J^{a}_{\gamma a})(J^{a}_{\beta\gamma a})^{-1}=
J^a_{\beta a}J^{\beta a}_{\beta\gamma a}J^{\beta\gamma a}_a=1.\]
Second, $c^1_a$ and $c^1_b$ are homologous:
\[c^1_a(\gamma)=J^a_{\gamma a}=J^a_bJ^b_{\gamma b}J^{\gamma b}_{\gamma a}=
J^a_b c^1_b(\gamma)(\gamma\cdot J^a_b)^{-1}.\text{ \qed}\]

\begin{definition}
 A non-commutative modular symbol is a non-abelian cohomology class in $H^1(\Gamma,\Pi)$, with representative \[c^1_a(\gamma)=J^a_{\gamma a},\]
\end{definition}

\section{Iterated integrals on membranes}

Iterated integrals on membranes are a higher dimensional analogue of iterated path integrals. 
This technical tool was used in \cite{MDZF} for constructing multiple Dedekind zeta values and in \cite{rec2} for proving new and classical reciprocity laws on algebraic surfaces.
It appeared first in the author preprint \cite{ModSym} for the purpose of non-commutative Hilbert modular symbols.

\subsection{Definition and properties}
Let $\H^1$ be the upper half plane. Let $\H^2$ be a product of two upper half planes. We are interested in the action of $GL_2(K)$, where $K$ is a real quadratic field. This group acts on $\H^2$ by linear fractional transforms. It is convenient to introduce cusp points $\P^1(K)$ as boundary points of $\H^2$.

Let $\omega_1,\dots,\omega_m$ be holomorphic $2$-forms on $\H^2$, which are continuous at the 
cusps $\P^1(K)$.
Let \[g:[0,1]^2\rightarrow \H^2\cup \P^1(K)\] be a continuous map, which is smooth almost everywhere. Denote by $F^1$ and $F^2$ the following coordinate-wise foliations: For any $a\in[0,1]$, define the leaves
\[F^1_a=\{(t_1,t_2)\in [0,1]^2\,|\,\,t_1=a\,\}.\]
and
\[F^2_a=\{(t_1,t_2)\in [0,1]^2\,| \,\,t_2=a\,\}.\]
\begin{definition}
We call the above map $g:[0,1]^2\rightarrow \H^2\cup\P^1(K)$ a membrane on $\H^2$ 
 if it is continuous and piecewise differentiable map such that $g(F^1_a)$ and $g(F^2_a)$ belong to a finite union of holomorphic curves in $\H^2\cup \P^1(K)$ for all constants $a$.
\end{definition}

Similarly, we define a membrane of a Hilbert modular variety.
Let $\omega_1,\dots,\omega_m$ be holomorphic $2$-forms on $Y_\Gamma=\H^2/\Gamma$, which are continuous at the cusps $\P^1(K)/\Gamma$. 
Let \[g:[0,1]^2\rightarrow X_\Gamma\]
be a continuous map, which is smooth almost everywhere,\ where $X_\Gamma=\H^2\cup \P^1(K)$ 
Let $f_i:X_\Gamma\rightarrow \P^1(\C)$ for $i=1,2$ be two algebraically independent rational functions on the Hilbert modular surface $X_\Gamma$.
Denote by $F^1$ and $F^2$ the following coordinate-wise foliations: For any $a\in[0,1]$, define the leaves
\[F^1_a=\{(t_1,t_2)\in [0,1]^2\,|\,\,t_1=a\,\}.\]
and
\[F^2_a=\{(t_1,t_2)\in [0,1]^2\,| \,\,t_2=a\,\}.\]
Let also
\[P^1_x=\{P\in X_\Gamma\,|\,f_1(P)=x\}.\]
and
\[P^2_x=\{P\in X_\Gamma\,|\,f_2(P)=x\}.\]
\begin{definition}
We call the above map $g:[0,1]^2\rightarrow X_\Gamma$  on $X_\Gamma$ if it is continuous and piecewise differentiable map such that for each $a$ there are $x_1$ and $x_2$ such that 
$g(F^1_a)\subset P^1_{x_1}$ and $g(F^2_a)\subset P^2_{x_2}$.
\end{definition}

We define three types of iterated integrals over membranes - type {\bf{a}}, type {\bf{b}} and type {\bf{c}}.
Type {\bf{a}} consists of linear iterations and type {\bf{b}} is more general and involves permutations. 
Type {\bf{a}} is less general, but more intuitive. The advantage of type {\bf{b}} is that it satisfies integral shuffle relation (Theorem \ref{thm shuffle}). In other words a product of two integrals of type {\bf{b}} can be expresses as a finite sum of iterated integrals over membranes of type {\bf{b}}. However, one might not be able to express a product of two integrals of type {\bf{a}} as a sum of finitely many integrals of type {\bf{a}}. Both type {\bf{a}} and type {\bf{b}} are defined on a product of two upper half planes. Type {\bf{c}} is defined on a Hilbert modular surface; that is, on a quotient of a product of upper half planes by an arithmetic group, which is commensurable to $SL_2({\cal{O}}_K)$. Type {\bf{c}} also satisfies a shuffle product, that is a product of two integrals of this type can be expresses a finite sum of such.

\begin{definition} 
\label{def membranes 2a}
 (Type {\bf{a}}, ordered iteration over membranes)
 Let \[g:[0,1]^2\rightarrow \H^2\cup \P^1(K)\]  be a membrane on $\H^2\cup \P^1(K)$.
\[\int_g\omega_1\cdots\omega_m
=
\int_{D}\bigwedge_{j=1}^m g^*\omega_i(t_{1,j},t_{2,j}),\]
where
\[
D=\{(t_{1,1},\dots, t_{2,m})\in[0,1]^{2m}\,\,\,|\,\,\,0\leq t_{1,1}\leq \dots\leq t_{1,m}\leq 1,\,\,\,
0\leq t_{2,1}\leq \dots\leq t_{2,m}\leq 1\}.
\]
\end{definition}

\begin{definition}
\label{def membranes 1b}
(Type {\bf{b}}, $2$ permutations)
 Let \[g:[0,1]^2\rightarrow \H^2\cup \P^1(K)\]  be a membrane on $\H^2\cup \P^1(K)$.
Let $\rho_1,\rho_2$ be two permutations of the set $\{1,2\dots,m\}.$
\[\int^{\rho_1,\rho_2}_g\omega_1\cdots\omega_m
=
\int_{D}\bigwedge_{j=1}^m g^*\omega_j(t_{1,\rho_1(j)},t_{2,\rho_2(j)}),\]
where
\[
D=\{(t_{1,1},\dots, t_{2,m})\in[0,1]^{2m}\,\,\,|\,\,\,0\leq t_{1,1}\leq \dots\leq t_{1,m}\leq 1,\,\,\,
0\leq t_{2,1}\leq \dots\leq t_{2,m}\leq 1\}.
\]

\end{definition}

\begin{definition}
\label{def membranes 1c}
(Type {\bf{c}}, $2$ permutations)
 Let \[g:[0,1]^2\rightarrow X_\Gamma\]  be a membrane on the Hilbert modular surface 
 $X_\Gamma=(\H^2\cup \P^1(K))/\Gamma$.
Let $\rho_1,\rho_2$ be two permutations of the set $\{1,2\dots,m\}.$
\[\int^{\rho_1,\rho_2}_g\omega_1\cdots\omega_m
=
\int_{D}\bigwedge_{j=1}^m g^*\omega_j(t_{1,\rho_1(j)},t_{2,\rho_2(j)}),\]
where
\[
D=\{(t_{1,1},\dots, t_{2,m})\in[0,1]^{2m}\,\,\,|\,\,\,0\leq t_{1,1}\leq \dots\leq t_{1,m}\leq 1,\,\,\,
0\leq t_{2,1}\leq \dots\leq t_{2,m}\leq 1\}.
\]
\end{definition}

{\bf{Examples of iterated integral of type b}}:
Let $\alpha_i(t_1,t_2)=g^*\omega_i(t_1,t_2)$. Denote by $(1)$ the trivial permutation and by $(12)$ the permutation exchanging $1$ and $2$.

1. The following $4$ diagrams 
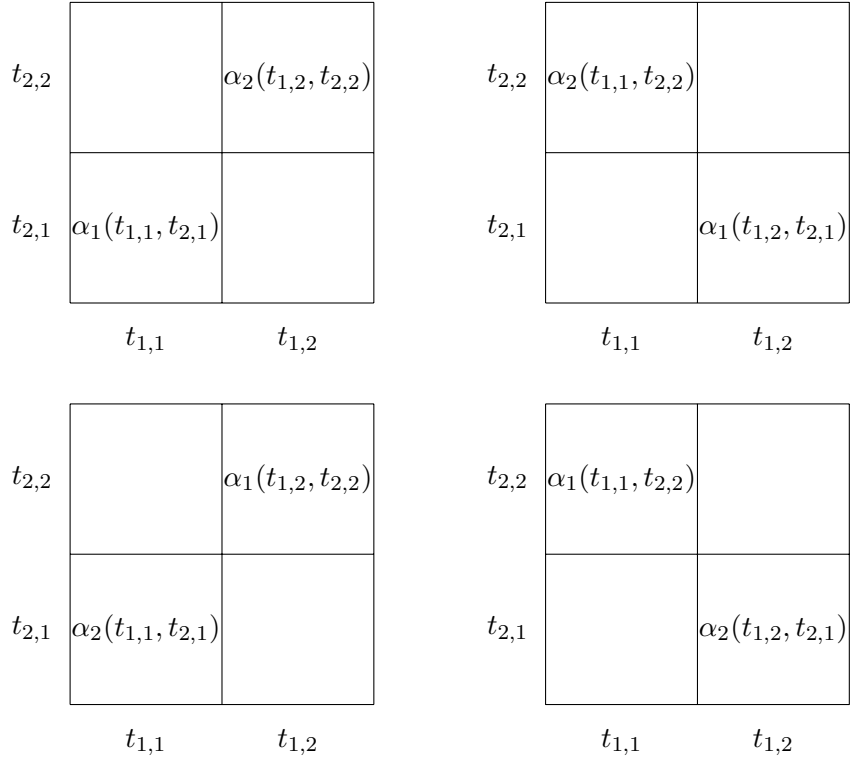
\vspace{.5cm}

%\begin{center}
\begin{tikzpicture}
\draw[step=2cm] (0,0) grid (4,4);
\draw (1,-.5)node{$t_{1,1}$};
\draw (-0.5,1)node{$t_{2,1}$};
\draw (3,-.5)node{$t_{1,2}$};
\draw (-0.5,3)node{$t_{2,2}$};
\draw (1,1)node{$\alpha_1(t_{1,1},t_{2,1})$};
\draw (3,3)node{$\alpha_2(t_{1,2},t_{2,2})$};
\end{tikzpicture}
%\end{center}
\hspace{1cm}
\begin{tikzpicture}
\draw[step=2cm] (0,0) grid (4,4);
\draw (1,-.5)node{$t_{1,1}$};
\draw (-0.5,1)node{$t_{2,1}$};
\draw (3,-.5)node{$t_{1,2}$};
\draw (-0.5,3)node{$t_{2,2}$};
\draw (3,1)node{$\alpha_1(t_{1,2},t_{2,1})$};
\draw (1,3)node{$\alpha_2(t_{1,1},t_{2,2})$};
\end{tikzpicture}

\vspace{.5cm}
%\begin{center}
\begin{tikzpicture}
\draw[step=2cm] (0,0) grid (4,4);
\draw (1,-.5)node{$t_{1,1}$};
\draw (-0.5,1)node{$t_{2,1}$};
\draw (3,-.5)node{$t_{1,2}$};
\draw (-0.5,3)node{$t_{2,2}$};
\draw (3,3)node{$\alpha_1(t_{1,2},t_{2,2})$};
\draw (1,1)node{$\alpha_2(t_{1,1},t_{2,1})$};
\end{tikzpicture}
%\end{center}
\hspace{1cm}
\begin{tikzpicture}
\draw[step=2cm] (0,0) grid (4,4);
\draw (1,-.5)node{$t_{1,1}$};
\draw (-0.5,1)node{$t_{2,1}$};
\draw (3,-.5)node{$t_{1,2}$};
\draw (-0.5,3)node{$t_{2,2}$};
\draw (1,3)node{$\alpha_1(t_{1,1},t_{2,2})$};
\draw (3,1)node{$\alpha_2(t_{1,2},t_{2,1})$};
\end{tikzpicture}

correspond, respectively, to the integrals 
\[\int_g^{(1),(1)}\omega_1\cdot\omega_2,  \,\,\,\,\,\,\,\, \int_g^{(12),(1)}\omega_1\cdot\omega_2,\]
\[\int_g^{(12),(12)}\omega_1\cdot\omega_2,\,\,\,\,\,\,\,\, \int_g^{(1),(12)}\omega_1\cdot\omega_2,\]

2. The following diagram

\begin{center}
\begin{tikzpicture}
\draw[step=2cm] (0,0) grid (6,6);
\draw (1,-.5)node{$t_{1,1}$};
\draw (3,-.5)node{$t_{1,2}$};
\draw (5,-.5)node{$t_{1,3}$};
\draw (-0.5,1)node{$t_{2,1}$};
\draw (-0.5,3)node{$t_{2,2}$};
\draw (-0.5,5)node{$t_{2,3}$};
\draw (3,1)node{$\alpha_1(t_{1,2},t_{2,1})$};
\draw (1,3)node{$\alpha_2(t_{1,1},t_{2,2})$};
\draw (5,5)node{$\alpha_2(t_{1,3},t_{2,3})$};
\end{tikzpicture}
\end{center}
corresponds to the integral
\[\int_g^{(12),(1)}\omega_1\cdot\omega_2\cdot\omega_3.\]

\begin{remark}
\label{remark1}
Let us give more intuition about Definition \ref{def membranes 1b}. Each of the differential forms 
$g^*\omega_1,\dots,g^*\omega_m$ has two arguments. Consider the set of first arguments for each of the  differential forms $g^*\omega_1,\dots,g^*\omega_m$. They are ordered as follows
\begin{equation}
\label{eq order 1}
0<t_{1,1}<t_{1,2}<\dots<t_{1,m}<1,
\end{equation}
(They are the coordinates of the domain $D$.) 
Since $g^*\omega_j$ depends on $t_{1,\rho_1(j)}$, we have that $t_{1,k}$ is an argument of 
$g^*\omega_{\rho_1^{-1}(k)}$, where $k=\rho_1(j)$. Then we can order the differential forms
$g^*\omega_1,\dots,g^*\omega_m$ according to the order of their first arguments given by the Inequalities \eqref{eq order 1}, which is
\[g^*\omega_{\rho_1^{-1}(1)},g^*\omega_{\rho_1^{-1}(2)},\dots,g^*\omega_{\rho_1^{-1}(m)}
\]
Similarly, we can order the differential forms $g^*\omega_1,\dots,g^*\omega_m$, with  respect to the order of their second arguments
\[g^*\omega_{\rho_2^{-1}(1)},g^*\omega_{\rho_2^{-1}(2)},\dots,g^*\omega_{\rho_2^{-1}(m)}
\]
We call the first ordering {\bf{horizontal}} and the second ordering {\bf{vertical}}.
\end{remark}

Now we are going to examine homotopy of a domain of integration and how that reflects on the integral.
Let $g_s:[0,1]^2\rightarrow \H^2\cup\P^1(K)$ be a family of membranes such that $g_s(0,0)=\infty$ and $g_s(1,1)=0$. Assume that the parameter $s$  is in the interval $[0,1]$.

Put $h(s,t_1,t_2)=g_s(t_1,t_2)$ to be a homotopy between $g_0$ and $g_1$.
Let 
\[G_s:[0,1]^{2m}\rightarrow  (\H^2\cup\P^1(K))^m,\]
be the map
\[G_s(t_{1,1},\dots,t_{2,m})=\left(g_s\left(t_{1,\sigma_1(1)},t_{2,\sigma_2(1)}\right),
g_s\left(t_{1,\sigma_1(2)},t_{2,\sigma_2(2)}\right),\dots,g_s\left(t_{1,\sigma_1(m)},t_{2,\sigma(m)}\right)\right).\]

Let $H$ be the induced homotopy between $G_0$ and $G_1$, defined by 
\[H(s,t_{1,1},\dots,t_{2,m})=G_s(t_{1,1},\dots,t_{2,m}).\]

We define diagonals in the domain $D\subset (0,1)^{2m}$, where 
\begin{align*}
D=\{(t_{1,1},t_{2,1},\dots,t_{1,m},t_{2,m})\in(0,1)^{2m}\,|\,
 &0\leq t_{1,1}\leq t_{1,2}\leq \dots\leq t_{1,m}\leq 1,\\
&\mbox{and}\,  0\leq t_{2,1}\leq t_{2,2}\leq \dots\leq t_{2,m}\leq 1\}.
\end{align*}
We define $D_{1,k}$ for $k=0,\dots,m$ as $D_{1,0}=D|_{t_{1,1}=0}$, $D_{1,k}=D|_{t_{1,k}=t_{1,k+1}}$, for $k=1,\dots,m-1$ and $D_{1,m}=D|_{t_{1,m}=1}$.
Similarly, we define
$D_{2,k}$ for $k=0,\dots,m$ as $D_{2,0}=D|_{t_{2,1}=0}$, $D_{2,k}=D|_{t_{2,k}=t_{2,k+1}}$, for $k=1,\dots,m-1$ and $D_{2,m}=D|_{t_{2,m}=1}$.

For iterated integrals of types {\bf{a}} and {\bf{b}}, we define diagonals in $V=(\H^2\cup \P^1(K))^m$. We denote a generic coordinate of 
$V=(\H^2\cup \P^1(K))^m$ by $(z_{1,1},z_{2,1},\dots,z_{1,m},z_{2,m})$
For $k=1,\dots,m-1$, let $V_{1,k}=V|_{z_{1,k}=z_{1,k+1}}$. Let also,
$V_{1,0}=V|_{z_{1,1}=0}$ and $V_{1,m}=V|_{z_{1,m}=1}$.
Similarly, for $k=1,\dots,m-1$, let $V_{2,k}=V|_{z_{2,k}=z_{2,k+1}}$. Let also,
$V_{2,0}=V|_{z_{2,1}=0}$ and $V_{2,m}=V|_{z_{2,m}=1}$.

For iterated integrals of type {\bf{c}}, we define ``diagonals" as fibers product of schemes corresponding to certain varieties. (for fiber product of schemes one may look at the book \cite{Hartshorne}. Occasionally, it will be more natural to realize the multiple fiber products as a finite limit in the category of schemes of finite type over $\C$. Let $X_{i,j}=X_\Gamma$ for $i,j=1,\dots,n$.
Let $V$ be the universal scheme (finite limit) that maps to
$X_{ij}$ for each $i$ and $j$ as a part of a commutative diagram. The commutative diagram is defined as follows: $X_{i,j}$ and $X_{i+1,j}$ both map to $\P^1(\C)$ via the morphism $f_1$ for $1\leq i\leq n-1$ and all $j$
and
$X_{i,j}$ and $X_{i,j+1}$ both map to $\P^1(\C)$ via the morphism $f_2$ for $1\leq j\leq n-1$ and all $i$.
Let $V_{1,0}$ be the subscheme of $V$ defined by putting $\P^1(\C)$ in the place of $X_{1,j}$, so 
that $f_1:X_{1,j}\rightarrow \P^1(\C)$ is replaced by the identity map and and the corresponding 
$f_2:X_{1,j}\rightarrow \P^1(\C)$ is deleted.
Let $V_{2,0}$ be the subscheme of $V$ defined by putting $\P^1(\C)$ in the place of $X_{i,1}$, so 
that $f_2:X_{1,j}\rightarrow \P^1(\C)$ is replaced by the identity map and and the corresponding 
$f_1:X_{1,j}\rightarrow \P^1(\C)$  is deleted.
Let $V_{1,n}$ be the subscheme of $V$ defined by putting $\P^1(\C)$ in the place of $X_{n,j}$, so 
that $f_1:X_{n,j}\rightarrow \P^1(\C)$ is replaced by the identity map and and the corresponding 
$f_2:X_{n,j}\rightarrow \P^1(\C)$  is deleted.
Let $V_{2,n}$ be the subscheme of $V$ defined by putting $\P^1(\C)$ in the place of $X_{i,n}$, so 
that $f_2:X_{n,j}\rightarrow \P^1(\C)$ is replaced by the identity map and and the corresponding 
$f_1:X_{n,j}\rightarrow \P^1(\C)$  is deleted.
Let also $V_{1,i}$ be the subscheme of $V$ obtained by 
replacing each factor $X_{i,j}\times_{\P^1(\C)}X_{i+1,j}$
by the corresponding diagonal for fixed $i$ and for all $j$.
And finally, let $V_{2,j}$ be the subscheme of $V$ obtained 
by replacing each factor $X_{i,j}\times_{\P^1(\C)}X_{i,j+1}$
by the corresponding diagonal for fixed $j$ and all $i$.

\begin{theorem}(Homotopy Invariance Theorem I)
\label{thm h invariance}
The iterated integrals on membranes from Definition \ref{def membranes 1b} (of type {\bf{b}})
are homotopy invariant, when the homotopy preserves the boundary of the membrane.\end{theorem}
\proof
Let 
\[\Omega=
\bigwedge_{j=1}^m\omega_j(z_{1,\sigma_1(j)},z_{2,\sigma_2(j)}).\]
Note that $\Omega$ is a closed form, 
since $\omega_i$ is a form of top dimension.
By Stokes Theorem, we have
\begin{align}
\label{eq g1 g0}
 0&=\int_{s=0}^{s=1}\int_D H^*d\Omega=\\
\label{eq boundary0}
& =\int_{D}G^*_1\Omega-\int_{D}G^*_0\Omega\pm\\
\label{eq boundary1}
&\pm\int_{s=0}^{s=1}\sum_{k=1}^{m-1}\left(\int_{D_{1,k}}\pm\int_{D_{2,k}}\right)H^*\Omega\\
\label{eq boundary2}
&\pm\int_{s=0}^{s=1}\left(\int_{D_{1,0}}\pm\int_{D_{2,0}}\right)H^*\Omega\\
\label{eq boundary3}
&\pm\int_{s=0}^{s=1}\left(\int_{D_{1,m}}\pm\int_{D_{2,m}}\right)H^*\Omega
\end{align} 
We want to show that the difference in the terms in \eqref{eq boundary0} is zero. 
It is enough to show that each of the terms \eqref{eq boundary1}, \eqref{eq boundary2} 
and \eqref{eq boundary3}  are zero. 
If $z_{1,k}=z_{1,k+1}$ for types {\bf{a}} and {\bf{b}} (or on $V_{1,k}$ for type {\bf{c}}), then the wedge of the corresponding differential forms will vanish. 
Thus the terms in \eqref{eq boundary1} are zero. 
If $z_1=0$ then $dt_1=0$, defined via the pull-back $H^*$. 
Then  the terms \eqref{eq boundary2} are equal to zero. 
Similarly, we obtain that the last integral 
\eqref{eq boundary3} vanishes. \qed

Let $A$ be a manifold with corners of dimension $2$ in $[0,1]^2$.
We recall the domain of integration
\[
D=\{(t_{1,1},\dots,t_{2,m})\in [0,1]^{2m}\,\,\,|\,\,\,0\leq t_{1,1} \leq \dots \leq t_{1,m} \leq 1,\,\,\,
0\leq t_{2,1}\leq \dots\leq t_{2,m}\leq 1\}.
\]
Let us define
\[
A^D=\{(t_{1,1},\dots,t_{2,m})\in D\,\,| \,\,  (t_{1,i},t_{2,j})\in A\mbox{ for }i,j=1,\dots,m\,\}\]
Let $\rho_1$ and $\rho_2$ be two permutations of $m$ elements. We define
\[G(t_{1,1},\dots,t_{2,m})=\left(g\left(t_{1,\rho_1(1)},t_{2,\rho_2(1)}\right),
g\left(t_{1,\rho_1(2)},t_{2,\rho_2(2)}\right),\dots,g\left(t_{1,\rho_1(m)},t_{2,\rho(m)}\right)\right).\]
as a function on $A^D$.
Recall
\[\Omega=
\bigwedge_{j=1}^m\omega_j(z_{1,\rho_1(j)},z_{2,\rho_2(j)}).\]
 
\begin{definition}
With the above notation, 
we define an iterated integral over a membrane of type {\bf{b}} restricted to a domain $U$, where $U=g(A)$ as
\[^b\int_{g,U}^{\rho_1,\rho_2}\omega_1\cdots\omega_m=\int_{A^D}G^*\Omega.\]
\end{definition}

Now we are going to define iterated integrals of type {\bf{c}}. 
\begin{definition}
Let $\Omega_0=\bigwedge_{i,j=1}^m\Omega_{i,j}$, where $\Omega_{i,j}=\omega_i\delta_{i,j}$ on 
$X_{i,i}\equiv X_\Gamma$, and where
$\Omega_{i,j}=1$ for $i\neq j$.
Let $in:X\rightarrow \prod_{i,j=1}^nX_{i,j}$
be the inclusion of the finite limit into the product of the schemes $X_{i,j}$.
Let $\Omega=in^*\Omega_0$.

With this definition of $\Omega$, we define iterated integrals of type {\bf{c}} restricted to a domain $U$, where $U=g(A)$ as
\[^c\int_{g,U}^{\rho_1,\rho_2}\omega_1\cdots\omega_m=\int_{A^D}G^*\Omega.\]
\end{definition}

Let $A_1$ and $A_2$  be two manifolds with corners, with a common component of the boundary as subsets of $[0,1]^2$. Let $A=A_1\cup A_2$.
Let $s$ be a map of sets with values $1$ or $2$,
\[s:\{1,\dots,m\}\rightarrow \{1,2\}.\]
We define a certain set $A^{D}_s$ as a subset of $A^D$ in the following way:
Consider the image of the map $G$. It has $m$ coordinates. 
The first coordinate, $g\left(t_{1,\rho_1(1)},t_{2,\rho_2(1)}\right)$, will be restricted to the set $A_{s(1)}$. The second coordinate,
$g\left(t_{1,\rho_1(2)},t_{2,\rho_2(2)}\right)$,  will be restricted to
$A_{s(2)},\dots$ and the last $m$ coordinate $g\left(t_{1,\rho_1(m)},t_{2,\rho_2(m)}\right)$
will be restricted to $A_{s(m)}$. Formally, this can be written as
\[A^D_s=\{(t_{1,1},\dots,t_{2,m})\in A^D\,\,| \,\, (t_{1,\rho_1(i)},t_{2,\rho_2(i)})\in A_{s(i)}\mbox{ for }i=1,\dots,m\}.\]
Note that the image of the map $s$ is $1$ or $2$.
\begin{definition}
\label{def it int 2 domains} 
With the above notation, we define an iterated integral of type {\bf{b}} or {\bf{c}} over two domains $U_1$ and $U_2$, where $U_i=g(A_i)$ and $U=U_1\cup U_2$ by
\begin{equation}
\label{eq decomposition 2}
\int_{g,U,s}^{\rho_1\rho_2} \omega_1\cdots\omega_m
 = 
\int_{A^D_s}G^*\Omega.
\end{equation}\
For type {\bf{b}} we have that $U$ is in $\H^2\cup \P^1(K)$ and
for type {\bf{c}} we have that $U$ is in $X_\Gamma=(\H^2\cup \P^1(K))/\Gamma$.
\end{definition}

Again we examine homotopy of iterated integrals on membranes. Now we restrict the domain of integration to a manifold with corners $A$ as a subset of $[0,1]^2$.
Assume that for the boundary of a domain $A$, denoted by $\partial A$,  
we have that $g(\partial A)$ belongs to a finite union of complex analytic curves in $\H^2$ for type {\bf{b}} and in $X_\Gamma$ for type {\bf{c}}.  
We call a {\it{complex boundary}} of $g(\partial A)$ the minimal union of complex analytic (holomorphic) curves such that $g(\partial A)$ belongs to a finite union 
of complex analytic curves in $\H^2$ for type {\bf{b}} and in $X_\Gamma$ for type {\bf{c}}.

\begin{theorem} (Homotopy Invariance Theorem II)
\label{thm h-invariance II}
Iterated integrals over membranes are homotopy invariant with respect to a homotopy that changes the boundary $\partial U$ of the domain of integration $U$, so that the boundary varies on a finite union of  complex analytic curves.
\end{theorem}
\proof
Assume that $g_0(\partial A)$ and $g_1(\partial A)$ have the same complex boundary.
Let $h$ be a homotopy between $g_0$ and $g_1$, 
such that for each value of $s$ we have that $h(s,\partial A)$ has the same complex boundary  as $h(0,\partial A)=g_0(\partial A)$. 
Let $A\subset B$ be a strict inclusion of disks. 
Identify $B-A^\circ$ with $A\times [0,1]$. Let 
$i:B-A^\circ \rightarrow [0,1]\times \partial A.$
Here $A^\circ$ is the interior of $A$ and $\partial A$ is the boundary of $A$.
Let $\tilde{g}_0$ be a map from $B$ to $\H^2$ 
so that $\tilde{g}_0(a)=g_0(a)$ for $a\in A$ and $\tilde{g}_0(b)\in h(i(b))$.
Since the restriction of pull-back $(\tilde{g}_0^*\omega_i)|_{B-A}=0$ 
is mapped to a finite union of complex curves, then it vanishes. Therefore
\begin{equation}
\label{eq g0}
\int_Ag_0^*\Omega=\int_B\tilde{g}_0^*\Omega.
\end{equation}
Let $\tilde{g}_1$ be a membrane from $B$ defined by
$\tilde{g}_1(a)=g_1(a)$ for $a\in A$ and $\tilde{g}_1(b)=\tilde{g}_1(a)$ for i(b)=(s,a). (Note that $i(b)\in [0,1]\times \partial A$.)
Again
\begin{equation}
\label{eq g1}
\int_A g_1^*\Omega=\int_B\tilde{g}_1^*\Omega.
\end{equation}
However, the boundary of $B$ is mapped to the same set (point-wise) by both $\tilde{g}_0$ and $\tilde{g}_1$. Moreover, the homotopy between $g_0$ and $g_1$ extends to a homotopy between
$\tilde{g}_0$ and $\tilde{g}_1$ that respects the inclusion into the complex boundary. 
Thus by Theorem \ref{thm h invariance}, we have that
\[\int_B\tilde{g}_0^*\Omega=\int_B\tilde{g}_1^*\Omega.\]
Using Equations \eqref{eq g0} and \eqref{eq g1}, we complete the proof of this theorem. \qed

\subsection{Generating series}
We are going to define two types of generating series  - type {\bf{a}} and type {\bf{b}}, corresponding to the iterated integrals on membranes of type {\bf{a}} and type {\bf{b}}. 

\begin{definition} (Type a)
\label{def int type a}
Let $A$ be a domain in $\R^2$. Let $g$ be a membrane. Let $U=g(A)\subset \H^2$. 
And let $\omega_1,\dots,\omega_m$ be holomorphic 2-forms on $\H^2$.
We define a generating series of type {\bf{a}} by
\[J^a(U)=1+\sum_{k=1}^\infty\sum_{c:\{1,\dots,k\}\rightarrow\{1,\dots,m\}}
X_{c(1)}\otimes\dots\otimes X_{c(k)}\int_{g,U}\omega_{c(1)}\dots\omega_{c(k)},
\]
where $c:\{1,\dots,k\}\rightarrow\{1,\dots,m\}$ be a map of sets.
\end{definition}

Consider a map of sets $c:\{1,\dots,k\}\rightarrow\{1,\dots,m\}$ and two permutations $\rho_1$, $\rho_2$ of $\{1,2,\dots,k\}$. We call two triples $(c',\rho_1',\rho_2')$ and $(c'',\rho_1'',\rho_2'')$ equivalent if  they are in the same orbit of the permutation group $S_k$. That is, $(c'',\rho_1'',\rho_2'')\sim( c',\rho_1',\rho_2')$  if for some $\tau\in S_k$ we have
$c''=c'\tau^{-1}$, $\rho_1''=\rho_1'\tau^{-1}$ and $\rho_2''=\rho_2'\tau^{-1}$.
Then for each equivalence class of a triple $(c,\rho_1,\rho_2)$, we can associate a unique pair $(c\circ \rho_1, c\circ\rho_2)$, (which are precisely the indices of the $X$-variables and $Y$-variables in \eqref{eq int0} and \eqref{eq int}, respectively.) The reason for using such an equivalence is that the integral in \eqref{eq int} is invariant by the above action of $\tau\in S_k$ on the triple $(c,\rho_1,\rho_2)$. 

\begin{definition} (Ring R, values of the generating series)
\label{def ring R} The values of the generation series 
of iterated integrals on membranes will be in a ring $R$, which we define as follows. Let $R_0$
be the quotient of the ring of formal power series 
\[R_0=\C\ll X_1,Y_1,\dots,X_m,Y_m\gg/I\]
modulo the two-sided ideal $I$ generated by $X_iY_j-Y_jX_i$ for $i,j=1,\dots,m$.
Let $R\subset R_0$ be the subring of formal power series whose monomials have the following property: in every monomial of $R'$, $X_{i}$ occurs as many times as $Y_{i}$.
\end{definition}

\begin{definition}
\label{def it int on U type b}  
(Type b)
We define the generating series of type {\bf{b}} on $U$ by
\begin{align}
\label{eq int0}
J^b(U)
=&1+\sum_{k=1}^\infty\sum_{(c,\rho_1,\rho_2)/\sim}
X_{c(\rho_1^{-1}(1))}\otimes\dots\otimes X_{c(\rho_1^{-1}(k))}\otimes\\
\label{eq int}
&\otimes
Y_{c(\rho_2^{-1}(1))}\otimes\dots\otimes Y_{c(\rho_2^{-1}(k))}
\int^{\rho_1,\rho_2}_{g,U}\omega_{c(1)}\dots\omega_{c(k)},
\end{align}
where the second summation is over all maps of sets
$c:\{1,\dots,k\}\rightarrow\{1,\dots,m\}$ and all permutations $\rho_1$, $\rho_2$ of $k$ elements, 
up to the above equivalence.
\end{definition}

Let $Y_\Gamma$ be a Hilbert modular surface. Let $\alpha$ and $\beta$ be two rational functions on $Y_\Gamma$. We denote by $D$ the union of the divisors $(\alpha)_\infty$ and $(\beta)_\infty$ at infinity. Let $F:Y_\Gamma-D\rightarrow \C^2$ be defined as $F(y)=(\alpha(y),\beta(y))$.
Let $g:(0,1)^2\rightarrow Y_\Gamma-D$ be a membrane, so that the composition $F\circ g$ respects the coordinate-wise foliations. Consider the differential forms $\omega_i$ from the Definition of type b. They are invariant under the action of the arithmetic group $\Gamma$. Thus, we can treat them as differential forms on the Hilbert modular variety $Y_\Gamma$. 
\begin{definition}
\label{def it int on U type c} 
(Type c)
With the new definition of a membrane $g$, and a domain $U\subset Y_\Gamma$, we define 
the generating series of type {\bf{c}} by
\begin{align}
\label{eq int0}
J^c(U)
=&1+\sum_{k=1}^\infty\sum_{(c,\rho_1,\rho_2)/\sim}
X_{c(\rho_1^{-1}(1))}\otimes\dots\otimes X_{c(\rho_1^{-1}(k))}\otimes\\
\label{eq int}
&\otimes
Y_{c(\rho_2^{-1}(1))}\otimes\dots\otimes Y_{c(\rho_2^{-1}(k))}
\int^{\rho_1,\rho_2}_{g,U}\omega_{c(1)}\dots\omega_{c(k)},
\end{align}
where the second summation is over all maps of sets
$c:\{1,\dots,k\}\rightarrow\{1,\dots,m\}$ and all permutations $\rho_1$, $\rho_2$ of $k$ elements, 
up to the above equivalence.
\end{definition}
\begin{definition} (Ring $R'$, genrating series $J(U_1,U_2)$) 
We define a generating series of iterated integrals on two disjoint domain $U_1$ and $U_2$ 
(see Definition \ref{def it int 2 domains}). Let $U_i=g(A_i)$.
\begin{align}J(U_1,U_2)
=&
1
+
\sum_{k=1}^\infty
\sum_{s:\{1,\dots,k\}\rightarrow\{1,2\}}
\sum_{(c,\rho_1,\rho_2)/\sim}X_{c(\rho_1^{-1}(1)),s(1)}\otimes\dots\otimes X_{c(\rho_1^{-1}(k)),s(k)}\otimes\\
\label{eq int2}
&\otimes
Y_{c(\rho_2^{-1}(1)),s(1)}\otimes\dots\otimes Y_{c(\rho_2^{-1}(k)),s(k)}
\int^{\rho_1,\rho_2}_{g,U,s}\omega_{c(1)}\dots\omega_{c(k)},
\end{align}
The generating series takes values in a ring $R'$ defined as follows. Let $R'_0$ be a quotient of the ring of formal power series
\[R'_0=\C\ll X_{1,1},X_{1,2},Y_{1,1},Y_{1,2},\dots,X_{m,1},X_{m,2},Y_{m,1},Y_{m,2}\gg/I',\]
where $I'$ is the two-sided ideal generated by Lie commutators of $X$ with any subscript and $Y$ with any subscript. Let $R'$ be a subring of $R'_0$ with the property: in every monomial of $R'$, $X_{i,j}$ occurs as many times as $Y_{i,j}$.
\end{definition}

\begin{lemma}
\label{lemma R' R}
Let $\phi:R'\rightarrow R$ be a homomorphism of rings defined by
$\phi(X_{i,1})=\phi(X_{i,2})=X_i$ and $\phi(Y_{i,1})=\phi(Y_{i,2})=Y_i$.
If $U=U_1\cup U_2$ is in $\H^2\cup\P^1(K)$ then
\[\phi(J(U_1,U_2))=J^b(U).\]
If $U=U_1\cup U_2$ is in $X_\Gamma$ then
\[\phi(J(U_1,U_2))=J^c(U).\]
\end{lemma}
\proof After applying the homomorphism $\phi$ the formal variables on the left hand side become independent of the map $s$. Therefore, we have to examine what happens when we sum over all possible maps $s$. The map $s(i)$ is $1$ or $2$. Their meaning is the following: If $s(i)=1$ then we restrict the form $g^*\omega_{c(i)}$ to $A_1$ (instead of to $A$). Similarly, if $s(i)=2$, we restrict $g^*\omega_{c(i)}$ to $A_2$. If we add both choices, restriction to $A_1$ and restriction to $A_2$, then we obtain restriction of $g^*\omega_{c(i)}$ to $A=A_1\cup A_2$. Thus, we obtain the formula
\[
\sum_{s:\{1,\dots,k\}\rightarrow\{1,2\}}
\int^{\rho_1,\rho_2}_{g,U,s}\omega_{c(1)}\dots\omega_{c(k)}
=
\int^{\rho_1,\rho_2}_{g,U}\omega_{c(1)}\dots\omega_{c(k)}.
\]
We do the same for every monomial in $R$. That proves the above Lemma for the generating series. \qed

\subsection{Shuffle product of generating series}
The regions of integration that we are mostly interested in will be ideal diangles, that is, a $2$-cell whose boundary has two vertices and two edges, and ideal triangles. 
All other regions that we will deal with are going to be a finite union of ideal diangles and ideal triangles. 
The first type of decomposition is based on a union of two diangles with a common vertex. 
The second type of decomposition will be based on two of the cells (diangles, or triangles) with a common edge.

Let $g_1$ and $g_2$ be two membranes.
Let $P=(0,0)$ and $Q=(1,1)$ be the vertices of a diangle $A$ 
and $Q=(1,1)$ and $R=(2,2)$ are the two points of a diangle $B$ 
as subsets of $\R^2$.  Assume that 
$A$ lies within the rectangle with vertices $(0,0), (0,1), (1,1), (1,0)$. Similarly,
assume that
$B$ lies within the rectangle $(1,1),(1,2),(2,2), (2,1)$.
Let $U=g(A)$ and $V=g(B)$. 
\begin{theorem}
(i)\[\int_{g,U\cup V}\omega_1\dots\omega_m
=
\sum_{j=0}^m\int_{g,U}\omega_1\dots\omega_j
\int_{g,V}\omega_{j+1}\dots\omega_m;\]

(ii) The generating series of type {\bf{a}} from Definition \ref{def int type a} satisfies the following property:
\[J^a(g;A\cup B;\Omega)=J^a(g;A;\Omega)J^a(g;B;\Omega).\]
\end{theorem}
The proof of the first statement is essentially the same as the combinatorial proof for composition of path, when one considers iterated path integrals (see Lemma \ref{lemma product of paths}). 
The second statement is combining all compositions into generating series (see Definition \ref{def int type a}), resembling the Manin's approach for non-commutative modular symbol.

For generating series of type {\bf{b}}, we have a similar statement. 
%However, first, we need to define a co-product on iterated integrals over membranes of type b. Let $\rho$ be a permutation of $\{1,2\dots,m\}$. Then
%\begin{align}
%&\Delta(X_1\dots X_m\otimes Y_{\rho(1)}\dots Y_{\rho(m)}
%\int_{g,A}^\rho\omega_1\dots,\omega_m)=\\
%&=\sum'_j 
%(X_1\dots X_j)(X_{j+1}\dots X_m)\otimes \\
%&\otimes(Y_{\rho(1)}\dots Y_{\rho(j)})(Y_{\rho(j+1)}\dots Y_{\rho(m)})
%\int_{g,A}^\rho\omega_1\dots,\omega_j
%\int_{g,A}^\rho\omega_{j+1}\dots,\omega_m,
%\end{align}
%where the summation is over all indices $j$ such that $\rho$ restricted to the first $j$ elements 
%$\{1,2,\dots,j\}$ of the set $\{1,2,\dots,m\}$ is a permutation of the first $j$ elements.

%With respect to the co-product, we have

\begin{definition}
\label{def shuffle of permutations}
Let $\rho'$ and $\rho''$ be two permutations of the sets $\{1,\dots,i\}$, $\{i+1,\dots,i+j\}$, respectively.
We define the permutation $\rho'^{-1}\cup\rho''^{-1}$ of $\{1,\dots,i+j\}$, which acts on $\{1,\dots,i\}$ as $\rho'^{-1}$ 
and on $\{i+1,\dots,i+j\}$ as $\rho''^{-1}$. We define the set of shuffles of two given permutations, denoted 
by $sh(\rho',\rho'')$, as  the set of all permutations $\rho$ of the set $\{1,2,\dots,i+j\}$ such that $\rho^{-1}$ 
is the composition of a shuffle of sets $\tau\in sh(i,j)$ (see Definition \ref{def shuffle of sets})
and  with $\rho'^{-1}\cup\rho''^{-1}$. That is,
\[\rho^{-1}=\tau\circ(\rho'^{-1}\cup\rho''^{-1}).\]
\end{definition}

\begin{definition}
\label{def shuffle of monomials}
We define a shuffle of two monomials
\[M'
=
X_{c'(\rho_1'^{-1}(1))}\otimes\dots\otimes X_{c'(\rho_1'^{-1}(i))}\otimes 
Y_{c'(\rho_2'^{-1}(1))}\otimes\dots\otimes Y_{c'(\rho_2'^{-1}(i))}
\int^{\rho'_1,\rho'_2}_{g,U'}\omega_{c'(1)}\dots\omega_{c'(i)}
\]
and
\[M''
=
X_{c''(\rho_1''^{-1}(1))}\otimes\dots\otimes X_{c''(\rho_1''^{-1}(j))}\otimes 
Y_{c''(\rho_2''^{-1}(1))}\otimes\dots\otimes Y_{c''(\rho_2''^{-1}(j))}
\int^{\rho''_1,\rho''_2}_{g,U''}\omega_{c''(i+1)}\dots\omega_{c''(i+j)},
\]
where $\rho'_1$ and $\rho'_2$ are permutations of $\{1,\dots,i\}$ and $c'$ is a map of sets
$c':\{1,\dots,i\}\rightarrow \{1,\dots,m\}$, 
and
$\rho''_1$ and $\rho''_2$ are permutations of $\{i+1,\dots,i+j\}$ and $c''$ is a map of sets
$c'':\{i+1,\dots,i+j\}\rightarrow \{1,\dots,m\}$. 
By a shuffle product of the monomials $M'$ and $M''$, we mean the following sum
\begin{align*}
M'\times_{Sh}M''
=&\sum_{\rho_1\in sh(\rho_1',\rho_1''),\rho_2\in sh(\rho_2',\rho_2'')}
X_{c(\rho_1^{-1}(1)),s(1)}\otimes\dots\otimes X_{c(\rho_1^{-1}(i+j)),s(i+j)}\otimes \\
&\otimes Y_{c(\rho_2^{-1}(1)),s(1)}\otimes\dots\otimes Y_{c(\rho_2^{-1}(i+j)),s(i+j)}
\int^{\rho_1,\rho_2}_{g,U'\cup U'',s}\omega_{c(1)}\dots\omega_{c(i+j)},
\end{align*}
where $c:\{1,\dots,i+j\}\rightarrow \{1,\dots,m\}$ such that the map $c$ restricted to the first $i$ elements is $c'$ and $c$ restricted to the last $j$ elements is $c''$. Here the maps $s$ takes the value $1$ on the set $c^{-1}\{1,\dots,i\}=c'^{-1}\{1,\dots,i\}$ and it takes the value $2$ on the set $c^{-1}\{i+1,\dots,i+j\}=c''^{-1}\{i+1,\dots,i+j\}.$ 
\end{definition}

\begin{theorem} 
\label{thm shuffle}
(Shuffle product) For iterated integrals of type {\bf{b}} and the corresponding generating series, we have the following shuffle relations: 

(i)\begin{align}
\label{eq shuffle1}
\int^{\rho_1',\rho_2'}_{g,U}\omega_1\dots \omega_j\int^{\rho_1'',\rho_2''}_{g,U}\omega_{j+1}\dots\omega_m
%=&\sum_{\rho_1\in sh(\rho_1',\rho_1''),\rho_2\in sh(\rho_2',\rho_2'')}\\
%&\int_{g,A} ^{\rho_1,\rho_2}\left(\bigwedge_{i=1}^m\omega_{i}(t_{1,\rho_1(i)},t_{2,\rho_2(i)})\right)=\\
=&\sum_{\rho_1\in sh(\rho_1',\rho_1''),\rho_2\in sh(\rho_2',\rho_2'')}
\int_{g,U}^{\rho_1,\rho_2}\omega_1\cdots\omega_{m}
\end{align}

(ii)\begin{align}
\int^{\rho_1',\rho_2'}_{g,U'}\omega_1\dots \omega_j\int^{\rho_1'',\rho_2''}_{g,U''}\omega_{j+1}\dots\omega_m
=&\sum_{\rho_1\in sh(\rho_1',\rho_1''),\rho_2\in sh(\rho_2',\rho_2'')}
\int_{g,U,s}^{\rho_1,\rho_2}\omega_1\cdots\omega_{m},
%&=
%\sum_{\rho_1\in sh(\rho_1',\rho_1''),\rho_2\in sh(\rho_2',\rho_2'')}\\
%&\int_{g,A'\cup A''} ^{\rho_1,\rho_2}
%\left(\bigwedge_{i=1}^j\omega_{i,A'}(s_{\rho_1(i)},t_{\rho_2(i)})\right)
%\wedge\\
%&
%\wedge\left(\bigwedge_{i=j+1}^m\omega_{i,A''}(s_{\rho_1(i)},t_{\rho_2(i)})\right).
\end{align}
where $s$ is a map from $\{1,\dots,m\}$ to $\{1,2\}$ so that $\{1,\dots,j\}$ are mapped to $1$ and the remaining elements are mapped to $2$.

(iii) \begin{equation}
\phi(J^b(U')\times_{Sh}J^b(U''))=J^b(U'\cup U'')
\end{equation}

(iv) \begin{equation}
\phi(J^c(U')\times_{Sh}J^c(U''))=J^c(U'\cup U'')
\end{equation}
\end{theorem}
\proof
For part (i), it is useful to consider the two orders of differential forms, given in Remark \ref{remark1}.
Note that we need to order the forms both horizontally and vertically in the terminology of Remark 
\ref{remark1}. Let us consider first the horizontal order. That is the order with respect to the first 
variables of the differential forms $g^*\omega_{\rho_1'^{-1}(1)},\dots,g^*\omega_{\rho_1'^{-1}(j)}$ and
$g^*\omega_{\rho_1''^{-1}(j+1)},\dots,g^*\omega_{\rho_1''^{-1}(m)}$ , corresponding to the  two integrals on the left hand side of Equation \eqref{eq shuffle1}. In order to arrange both of the above orderings in one sequence of increasing first 
arguments, we need to shuffle them (similarly to a shuffle of a deck of cards.) 
That leads to  $\rho_1\in sh(\rho'_1,\rho''_1)$ (see Definition \ref{def shuffle of permutations}). 
We proceed similarly, with the second arguments and the permutations 
$\rho'_2,\rho''_2$ and $\rho_2$.

For part (ii) apply the equality from part (i) when the differential forms $g^*\omega_1,\dots,g^*\omega_j$ are multiplied by 
the function ${\bf{1}_{A'}}$ defined by 
\[{\bf{1}_{A'}}(x)
=
\left\{
\begin{tabular}{ll}
$1$ & for $x\in A'$\\
$0$ & for $x\notin A'$
\end{tabular}
\right.
\]
and the differential forms $g^*\omega_{j+1},\dots,g^*\omega_m$ are multiplied by ${\bf{1}_{A''}}$.
%Using  part (ii) together with Equation \eqref{eq decomposition 2}, we prove (iii).
For part (iii), we are going to establish similar relation among generating series as elements of $R'$. Applying the homomorphism $\phi:R'\rightarrow R$ from Lemma \ref{lemma R' R}, we obtain desired equality. Every monomial from $J(U_1)$ is of the form
\[M'
=
X_{c'(\rho_1'^{-1}(1))}\otimes\dots\otimes X_{c'(\rho_1'^{-1}(i))}\otimes 
Y_{c'(\rho_2'^{-1}(1))}\otimes\dots\otimes Y_{c'(\rho_2'^{-1}(i))}
\int^{\rho'_1,\rho'_2}_{g,U'}\omega_{c'(1)}\dots\omega_{c'(i)}
\]
and similarly every monomial from $J(U_2)$ is of the form
\begin{align*}
M''
=&
X_{c''(\rho_1''^{-1}(1))}\otimes\dots\otimes X_{c''(\rho_1''^{-1}(j))}\otimes \\
&\otimes Y_{c''(\rho_2''^{-1}(1))}\otimes\dots\otimes Y_{c''(\rho_2''^{-1}(j))}
\int^{\rho''_1,\rho''_2}_{g,U''}\omega_{c''(i+1)}\dots\omega_{c''(i+j)},
\end{align*}
where $\rho'_1$ and $\rho'_2$ are permutations of $\{1,\dots,i\}$ and $c'$ is a map of sets
$c':\{1,\dots,i\}\rightarrow \{1,\dots,m\}$, 
and
$\rho''_1$ and $\rho''_2$ are permutations of $\{i+1,\dots,i+j\}$ and $c''$ is a map of sets
$c'':\{i+1,\dots,i+j\}\rightarrow \{1,\dots,m\}$. 
We take the shuffle product of the monomials $M'$ and $M''$ (see Definition \ref{def shuffle of monomials})
\begin{align*}
M'\times_{Sh}M''
=&\sum_{\rho_1\in sh(\rho_1',\rho_1''),\rho_2\in sh(\rho_2',\rho_2'')}
X_{c(\rho_1^{-1}(1)),s(1)}\otimes\dots\otimes X_{c(\rho_1^{-1}(i+j)),s(i+j)}\otimes \\
&\otimes Y_{c(\rho_2^{-1}(1)),s(1)}\otimes\dots\otimes Y_{c(\rho_2^{-1}(i+j)),s(i+j)}
\int^{\rho_1,\rho_2}_{g,U,s}\omega_{c(1)}\dots\omega_{c(i+j)},
\end{align*}
where the map $s$ takes the value $1$ on the set $c^{-1}\{1,\dots,i\}$ and takes the value $2$ on the set $c^{-1}\{i+1,\dots,i+j\}.$ It determines the map $s$ uniquely.

In order to complete the proof, we have to show that every monomial in $J(U_1,U_2)$ can be obtained 
in exactly one way as a result (on the right hand side) 
of a shuffle product of a pair of monomials $(M_1,M_2)$ from $J(U_1)$ and $J(U_2)$.
Every monomial from $J(U_1,U_2)$ is characterized by two permutation $\rho_1,\rho_2$, 
and two maps of sets $c:\{1,\dots,k\}\rightarrow \{1,\dots,m\}$ and $s:\{1,\dots,k\}\rightarrow \{1,2\}$.
Let $i$ be the number of elements in $s^{-1}(1)$ and $j$ be the number of elements in $s^{-1}(2)$. Then $i+j=k$. Then $i$ is the number of differential forms among $g^*\omega_{c(1)},\cdots, g^*\omega_{c(k)}$, which are restricted to the set $A_1$. The remaining $j$ differential forms are restricted to $A_2$. 
Also, every permutation $\rho_1$ can be written in an unique way as a composition of a shuffle 
$\tau_1\in sh(i,j)$ and two disjoint permutations $\rho'_1$ and $\rho''_1$ of $i$ and of $j$ elements,
respectively (see Definition \ref{def shuffle of permutations}). 
Similarly, $\rho_2$ can be written in a unique way as a product of a shuffle $\tau_2\in sh(i,j)$ 
and two disjoint permutation $\rho_2'$ and $\rho_2''$. The map of sets $c_1$ is defined as a restriction of the map $c$ to the image of $\rho'_1$. Similarly, the map $c_2$  is defined as a restriction of the map $c$ to the image of $\rho''_1$. 
Now we can define the monomials $M'$ and $M''$ in $J(U_1)$ and $J(U_2)$, 
based on the triples $\rho'_1,\rho'_2,c'$ and $\rho''_1,\rho''_2,c''$, respectively.
Such monomials are unique. One can show that the shuffle product of $M'$ and $M''$ contains the monomial in $J(U_1,U_2)$, that we started with, exactly once.
The proof of part (iii) is complete after applying Lemma \ref{lemma R' R}. \qed

\section{Hilbert modular symbols}
In this Section, we recall the Hilbert modular group and its action on the product of two upper half planes. Then we define commutative Hilbert module symbol, (Subsection 3.1) and its pairing with the cohomology of the Hilbert modular surface, (Subsection 3.2). In Subsections 3.3 and 3.4, we define the non-commutative Hilbert module symbols (Definition \ref{def c1 c2}) as a generating series of iterated integrals over membranes of type {\bf{b}}. We also examine relations among the non-commutative Hilbert modular symbols (Theorem \ref{thm NC relations}), which we interpret as cocycle conditions or as a difference by a coboundary (Theorem \ref{thm c1 c2}). In Subsection 3.5, we consider a two-category $C$ with a sheaf $J$ on $C$. Then the non-commutative Hilbert modular symbol a sheaf on a two-category. This is done in order to give a plausible approach to defining a suitable non-commutative cohomology set. In Subsection 3.6, we make explicit computations and compare them to computations for multiple Dedekind zeta values.

\subsection{Commutative Hilbert modular symbols}
In this Subsection, we define a commutative Hilbert modular symbol, using geodesics, geodesic triangles and geodesic diangles. Then, we prove certain relations among the commutative Hilbert modular symbols, which are generalized to relations among non-commutative Hilbert modular symbols (Subsection 3.4).

Let $K=\Q(\sqrt{d})$ be a real quadratic extension of $\Q$.
Then the ring of integers in $K$ is
\[{\cal{O}}_K
=
\left\{
\begin{tabular}{lll}
$\Z[\frac{1+\sqrt{d}}{2}]$ &for $d=1\mod 4$,\\
\\
$\Z[\sqrt{d}]$ &for  $d=2,3 \mod 4$.
\end{tabular}
\right.\]
Then  $\Gamma=SL_2({\cal{O}}_K)$ is called a Hilbert modular group. 
Let $\gamma\in \Gamma$. 
We recall the action of $\gamma$ on a product of two upper half planes $\H^2$. 
Let
\[\gamma=
\gamma_1=
\left(
\begin{tabular}{ll}
$a_1$ & $b_1$\\
$c_1$ & $d_1$
\end{tabular}
\right).
\]
Let $a_2,b_2,c_2,d_2$ be the Galois conjugate of $a_1,b_1,c_1,d_1$, respectively.
Let us define $\gamma_2$ by
\[
\gamma_2=
\left(
\begin{tabular}{ll}
$a_2$ & $b_2$\\
$c_2$ & $d_2$
\end{tabular}
\right).
\]
Let $z=(z_1,z_2)$ be any point of the product of two upper half planes $\H^2$.  

%We define \[\gamma z=(\gamma_1 z_1, \gamma_2 z_2),\]
%where \[\gamma_1 z_1=\frac{a_1z_1+b_1}{c_1z_1+d_1}\mbox{ and }
%\gamma_2 z_2=\frac{a_2z_2+b_2}{c_2z_2+d_2}\]
%are linear fractional transforms. 

For an element $\gamma\in GL_2(K)$, we define the following action:
If $\det\gamma$ is totally positive, 
that is  $\det \gamma_1>0$ and $\det \gamma_2>0$, 
then the action of $\gamma$ on $z=(z_1,z_2)\in\H^2$ 
is essentially the same as for $\gamma\in SL_2(K)$, 
namely, \[\gamma z=(\gamma_1 z_1,\gamma_2 z_2),\]
where
\[\gamma_1 z_1=\frac{a_1z_1+b_1}{c_1z_1+d_1}\,\,\mbox{ and }\,\,
\gamma_2 z_2=\frac{a_2z_2+b_2}{c_2z_2+d_2}
\]
are linear fractional transforms. 
If $\det \gamma$ is totally negative, that is, $\det \gamma_1<0$ and $\det \gamma_2<0$, then we define 
\[\gamma z= \left(-\frac{a_1\overline{z}_1+b_1}{c_1\overline{z}_1+d_1}, 
-\frac{a_2\overline{z}_2+b_2}{c_2\overline{z}_2+d_2}\right).\]
Similarly if one of the embeddings of $\det \gamma$ is positive and the other is negative, for example, 
$\det \gamma_1>0$ and $\det \gamma_2<0$, such as $\det \gamma=\sqrt{d}$, then
\[\gamma z= \left(\frac{a_1z_1+b_1}{c_1z_1+d_1}, 
-\frac{a_2\overline{z}_2+b_2}{c_2\overline{z}_2+d_2}\right).\]

We add cusp points $\P^1(K)$ to $\H^2$. Then the quotient 
$SL_2({\cal{O}}_K)\backslash (\P^1(K)\cup\H^2)$ is compact.

We are going to examine carefully geodesics joining the cusps $0$, $1$  and $\infty$.

Let $z_0,z_1,z_\infty$ be three distinct cusp points. 
There is a unique element $\gamma\in PGL_2(K)$ such that that send $z_0,z_1$ and $z_\infty$ to $0,1$ and $\infty$, respectively.

%Let $\gamma_0,\gamma_1,\gamma_\infty$ be four group elements of the Hilbert modular group. Let $z_i=\gamma_i(\infty)$ Then there exists an element of $\gamma\in GL_2(K)$ such that
%$\gamma(z_0)=0$, $\gamma(z_1)=1$ and $\gamma(z_\infty)=\infty$. Namely,
%\[\gamma=
%\left(
%\begin{tabular}{llll}
%$z_1-z_\infty$ &$-z_0(z_1-z_\infty)$\\
%$z_1-z_0$       &$-z_\infty(z_1-z_0)$
%\end{tabular}
%\right)\]
Let \[i:\H\rightarrow \H^2\]
\[i(x)=(x,x)\]
be the diagonal map and  $\Delta$ be its image.
Consider the Hirzebruch-Zagier divisor $X=\gamma^*\Delta$. 
It is an analytic curve that passes through the points $z_0$, $z_1$ and $z_\infty$. 
Then $X$ is a holomorphic curve in $\H^2$ if $\det \gamma$ 
is totally positive or totally negative.
If $\det \gamma$ is not totally positive or totally negative, 
then $X$ is a holomorphic curve in $\H^1\times \overline{\H}^1\cup \P^1(K)$, 
in other words it is anti-holomorphic curve in $\H^2$, such as 
$z_1=-\overline{z}_2$. 
Let $\Delta_X=\gamma^*\Delta$ 
be the pull-back of the geodesic triangle $\Delta$ 
between the points $0,1,\infty$ in the analytic curve $X$.

Given four points on the boundary in $\H^2\cup \P^1(K)$,
 we are tempted to consider them as vertices of a geodesic tetrahedron in 
 $\H^2\cup \P^1(K)$, whose faces are triangles of the type 
$\Delta_X$. However, there is one problem that we encounter: 
Two distinct cusps could be connected by two different  geodesics in $\H^2\cup \P^1(K)$. 
In particular, two triangles from the faces of the ``tetrahedron" might not have a common edge, 
but only two common vertices. 
Thus, we are led to consider a {\it{thickened tetrahedron}} with two types of faces on the boundary: 
the first type is an ideal triangle that we have just defined and the other type is an ideal diangle - a union of geodesics connecting two fixed points, 
which has the homotopy type of a disc with two vertices and two edges. 
The two edges of an ideal diangle in the boundary of a thickened tetrahedron 
correspond to the two geodesics connecting the same two cusps, 
where two geodesics belong to the geodesic triangles that have the two cusps in common.

Let us describe a diangle $D_{0,\infty;1,\alpha}$ whose two vertices are $0$ and $\infty$ and whose two sides are geodesics that belong to each of the ideal triangles $0,1,\infty$ and  $0,\alpha,\infty$.
The geodesic $l_0$ between the points $0$ and $\infty$ that lie on the geodesic triangle $0,1,\infty$ can be parametrized in the following way: $\{(it,it)\,|\,t\in \R,\, t\geq 0\}\subset Im(\H)\times Im(\H).$
Here by $Im(\H)$ we mean the imaginary part of the upper half plane.
The element $\gamma\in \Gamma$ that sends  $0,\alpha,\infty$ to  $0,1,\infty$ is
$\gamma
=
\left(
\begin{tabular}{ll}
$\alpha^{-1}$ & $0$\\
$0$ &$1$
\end{tabular}
\right)
$
Then $(\alpha^{-1})^*(it,it)=(|\alpha_1|it,|\alpha_2|it)$.
Therefore, the geodesic $l_\alpha$ between the points $0$ and $\infty$ that lie on the geodesic triangle $0,\alpha,\infty$ can be parametrized in the following way $\{(|\alpha_1|it,|\alpha_2|it)\,|\,t\in \R,\, t\geq 0\}\subset Im(\H)\times Im(\H).$ Then, we define the diangle $D_{0,\infty;1,\alpha}$ as the two dimensional region in $Im(\H)\times Im(\H)$ between the lines $l_0$ and $l_\alpha$. We also consider the diangle with orientation. If $|\alpha_1|>|\alpha_2|$ then it is positively oriented. If the inequality is reversed then the diangle is negatively oriented; if $|\alpha_1|=|\alpha_2|$  then it is a degenerate diagle, which consists of a single geodesic. 
All other diangles that we will consider are translates of $D_{0,\infty;1,\alpha}$  via the action of any element $\gamma\in PGL_2(K)$.

\begin{lemma}
(i) Each geodesic triangle  $\Delta_X$  lies either on a holomorphic curve or on an anti-holomorphic curve.

(ii) Each geodesic in a geodesic triangle $\Delta_X$ 
belongs both to a holomorphic curve and to an anti-holomorphic curve.
%(iii) Each geodesic diangle $D_T$ lies either in a holomorphic curve on in an anti-holomorphic curve.
\end{lemma}
Part (i) follows from the construction of a geodesic triangle before the lemma. 
For part (ii), consider the following:
Let $\Delta(0,1,\infty)$ be the geodesic triangle in the diagonal of $\H^2$ 
connecting the points $0$, $1$ and $\infty$. It is a holomorphic curve. 
Thus, a geodesic $\{(it,it)\in \H^2\,\,|\,\, t>0\}$, connecting the points $0$ and $\infty$ as a face of the 
geodesic triangle $\Delta(0,1,\infty)$ lies on a holomorphic curve. 
Now consider the geodesic triangle $D(0,\sqrt{d},\infty)$. It lies on an anti-holomorphic curve in 
$\H^2$, by which we mean a complex curve in $\H^2$, where we have taken the complex conjugate complex structure in one of the upper half planes. 
Since, the linear fractional transform that sends $D(0,\sqrt{d},\infty)$ to $D(0,1,\infty)$ 
does not have totally positive (or totally  negative) determinant. 
Explicitly, the linear fractional transform that sends $(0,\sqrt{d},\infty)$ to $(0,1,\infty)$ is 
\[\gamma=
\left(
\begin{tabular}{llll}
$1$ &$0$\\
$0$ &$\sqrt{d}$
\end{tabular}
\right)\]
Then 
\[(\gamma_1,\gamma_2)=
\left(\left(
\begin{tabular}{llll}
$1$ &$0$\\
$0$ &$\sqrt{d}$
\end{tabular}\right)
,
\left(
\begin{tabular}{llll}
$1$ &$0$\\
$0$ &$-\sqrt{d}$
\end{tabular}\right)\right)\]
We have 
$\gamma_1(it)=\frac{1}{\sqrt{d}}it$
and
$\gamma_2(it)=-\frac{1}{\sqrt{d}}\overline{it}=\gamma_1(it)$.
Then the same geodesic $(it,it)$ 
belongs to the anti-holomorphic curve given by the pull-back of the diagonal 
with respect to the linear fractional map $\gamma$. 
Thus, we obtain that the geodesic $(it,it)$, connecting $0$ and $\infty$, 
belongs to both a holomorphic curve and an anti-holomorphic curve. 
Similarly, any translate of the geodesic $(it,it)$ via a linear fractional map from $GL_2(K)$ 
would belong to both  a holomorphic curve and an anti-holomorphic curve. 
That proves part (ii). 
\qed

%For part (iii), we will examine the entire family of geodesics that connect $0$ and $\infty$.
%Consider the geodesic $(it,it)$ connecting $0$ and $\infty$. Consider those
%$A=(A_1,A_2)\in GL_2(\R)\times GL_2(\R)$ that fix both $0$ and $\infty$. 
%Then both $A_1$ and $A_2$ are diagonal. ???

\begin{definition}
Let
$p_1,p_2,p_3,p_4$ be cusp points in $\H^2\cup \P^1(K).$
To each triple of points $p_1,p_2,p_3$, we associate
the geodesic triangle $\{p_1,p_2,p_3\}$ with coefficient $1$ as an element of the singular chain complex in $C_2(\H^2\cup \P^1(K),\Q)$. Also, to each quadruple of points $p_1,p_2,p_3,p_4$,
we associate the geodesic diangle between the two geodesic connecting $p_1$ and $p_2$ so that the first geodesic is a face of the geodesic triangle $\{p_1,p_2,p_3\}$ and the second geodesic is a face of the geodesic triangle $\{p_1,p_2,p_4\}$. We denote such diangle by $\{p_1,p_2;p_3,p_4\}$.
We call the geodesic triangle $\{p_1,p_2,p_3\}$ and the geodesic diangle  $\{p_1,p_2;p_3,p_4\}$, 
considered as elements of $C_2(\H^2\cup \P^1(K),\Q)$, {\bf{commutative Hilbert modular symbols}}. 
\end{definition}

\begin{theorem}
The commutative Hilbert modular symbols modulo the boundary of singular $3$-chains, 
$\partial C_3(\H^2\cup \P^1(K),\Q)$ satisfy the following properties:

1. If $\sigma$ is a permutation of the set $\{1,2,3\}$ then  
\[\{p_{\sigma(1)},p_{\sigma(2)},p_{\sigma(3)}\} =sign(\sigma) 
\{p_1,p_2,p_3\}.\]

2. If $p_1,p_2,p_3,p_4$ are four points on the same holomorphic (or anti-holomorphic) curve  of the type $\gamma^*\Delta$ then
\[\{p_1,p_2,p_3\} + \{p_2,p_3,p_4\} = \{p_1,p_2,p_4\} + \{p_1,p_3,p_4\}.\]

To each four points  $p_1,p_2,p_3,p_4$, we associate a diangle with vertices $p_1$ and $p_2$. 
Let $\{p_1,p_2;p_3,p_4\}$ be the corresponding symbol.

3. If $p_1,p_2,p_3,p_4$ are four points on the same holomorphic (or anti-holomorphic) curve of the type $\gamma^*\Delta$ then
\[0=\{p_1,p_2;p_3,p_4\}.\]

4.  For every district four points $p_1,p_2,p_3,p_4$, we have the following relations based on the orientation of the domain
\[\{p_2,p_1;p_3,p_4\} =\{p_1,p_2;p_4,p_3\}  =-\{p_2,p_1;p_4,p_3\} = -\{p_1,p_2;p_3,p_4\}.\] 

5. For every five points $p_1,p_2,p_3,p_4,p_5$, we have 
 \[\{p_1,p_2;p_3,p_4\} + \{p_1,p_2;p_4,p_5\} = \{p_1,p_2;p_3,p_5\}.\]

6. We also have relation between the two types of commutative Hilbert modular symbols. 
For every four distinct points $p_1,p_2;p_3,p_4$, we have 
\begin{align*}
0=&\{p_1,p_2,p_3\} + \{p_2,p_3,p_4\} -\\
&- \{p_1,p_2,p_4\} - \{p_1,p_3,p_4\}+\\
&+\{p_1,p_2;p_3,p_4\} +\{p_2,p_3;p_1,p_4\} +\{p_3,p_1;p_2,p_4\}+\\
&+\{p_3,p_4;p_1,p_2\} +\{p_1,p_4;p_2,p_3\} +\{p_2,p_4;p_3,p_1\}.
\end{align*}
\end{theorem}
\proof Part 1 follows from orientation of the simplex in singular homology. Part 2 is an equality induced by two different triangulations on a holomorphic (or anti-holomoprhic) curve with $4$ vertices. In that setting the diangles are trivial, which proves Part 3. Part 4 follows from orientation of the diangle.
Part 5 corresponds to a union of two geodesic diangles with a common face, given by a third geodesic diangle. Part 5 will be used for a non-commutative $1$-cocycle relation for the non-commutative Hilbert modular symbol (see Conjecture \ref{thm c1 c2}). Part 6 is a boundary relation for the boundary of a thickened tetrahedron. By a thickened tetrahedron we mean a union of four geodesic triangles corresponding to each triple of points among the four points $p_1,p_2,p_3,p_4$ together with six  geodesic diangles that correspond to the area between the faces of the geodesic triangles. They correspond exactly to the thickening of the six edges of a tetrahedron. \qed

Part 6 will be used to derive explicit formulas  for the non-commutative Hilbert modular symbol of type {\bf{c'}} resembling  a non-commutative $2$-cocycle relation (see Conjecture \ref{thm c2}).

%\subsection{Homology of $\Gamma$}
%Let $\Gamma$ be a congruence subgroup of $SL_2({O}_K)$, where $K$ is a real quadratic field.
%The group acts on a product of two upper half planes together with the cusp points 
%$\H^2\cup\P^1(K)$. 
%We are going to consider the following geodesic ``triangulation" of $\H^2\cup\P^1(K)$. 
%The {\it{vertices}} will be all points of $\P^1(K)$, {\it{edges}} will be 
%all geodesics connecting two distinct cusp points in $\H^2\cup\P^1(K)$, 
%which are a face of an geodesic triangle described above, 
%{\it{faces}} will be all geodesic triangles and diangles describes above, 
%and finally $3-${\it{cells}} will be thickened tetrahedrons, 
%that is a tetrahedron, whose edges are thickened to become diangles, 
%which have two vertices, two edges and one face. 
%Let  $X=(\H^2\cup\P^1(K))/\Gamma$. 
%Let $C_*(\H^2\cup\P^1(K),\Q)$ be the chain of $\H^2\cup\P^1(K)$ 
%with respect to the above triangulation and  with rational coefficients. 
%Let $C_*(X,\Gamma)=H_0(C_*(\H^2\cup\P^1(K),\Q)$ 
%be the $\Gamma$-orbits of the elements of $C_*(\H^2\cup\P^1(K),\Q)$.
%\begin{lemma} The chain 
%$H_0(\Gamma,C_p(\H^2\cup\P^1(K),\Q))$ 
%computes the homology of $X$ up to degree $2$ with rational coefficients.
%\end{lemma}
%\proof
%Then by Hochschild-Serre spectral sequence, we have
%\[H_p(\Gamma,H_q(\H^2\cup\P^1(K),\Q)=>H_{p+q}(X,\Q).\]
%Then the homology of $X$ can be computed by the $\Gamma$-orbits 
%$H_0(\Gamma,C_p(\H^2\cup\P^1(K),\Q))$ for a degree less than or equal to $2$. 

\subsection{Pairing of the modular symbols with cohomology}

In this subsection, we consider pairings between commutative Hilbert modular symbols and cusp forms. In some cases, we prove that such pairings give periods in the sense of \cite{KZ}. 

We are interested in holomorphic cusp forms with respect to $\Gamma$. 
Equivalently, we can consider the holomorphic $2$-forms on $\tilde{X}$, 
which is the minimal smooth algebraic compactification of $X$ \cite{Hirz}. 
At this point we should distinguish between geodesic triangles $p_1,p_2,p_3$ 
that lie on a holomorphic curve or on anti-holomorphic curve. 
The reason for distinguishing is 
that a holomorphic $2$-form restricted to a holomorphic curve vanishes. 
The way to distinguish the two type of geodesic triangles is the following: 
Let $\gamma$ be a linear fractional transform that sends the points $p_1,p_2,p_3$ to $0,1,\infty$. 
If $\det \gamma$ is totally positive or totally negative 
then the geodesic triangle $p_1,p_2,p_3$ lies on a holomorphic curve. 
If $\det\gamma$ is not totally positive nor totally negative 
then the geodesic triangle $p_1,p_2,p_3$ lies on an anti-holomorphic curve. 

\begin{definition}
Let $M_2(\H^2\cup \P^1(K),\Q)$ be the span of the Hilbert modular symbols $\{p_1,p_2,p_3\}$ and $\{p_1,p_2;p_3,p_4\}$ as a subspace of the singular chain $C_2(\H^2\cup \P^1(K),\Q)$. 
We define the following pairing
\[<\,,\,>:M_2(\H^2\cup \P^1(K))\times S_{2,2}(\Gamma)\rightarrow \C,\]
by setting
\[<\{p_1,p_2,p_3\}, fdz_1\wedge dz_2>
= 
\int_{\{p_1,p_2,p_3\}}fdz_1\wedge dz_2\]
for  geodesic triangles and
\[<\{p_1,p_2;p_3,p_4\}, fdz_1\wedge dz_2>
=
\int_{\{p_1,p_2;p_3,p_4\}}fdz_1\wedge dz_2\]
for geodesic diangles. 
\end{definition}

We are going to use that a Hilbert modular surface $X(\C)$ can be realized as the complex points of an arithmetic surface defined over a number field $F$.
\begin{theorem}
\label{thm period}
The image of the above pairing is a period over a number field $F$, 
when we integrate a normalized cusp Hecke eigenform $f$ of weight $(2,2)$; (for Hecke eigenforms, see \cite{Shi}, \cite{Mladen}).
\end{theorem}
\proof From Lemma 3.1 (ii), the boundary of the geodesic triangles of the diangles are geodesics that lie on 
holomorphic curves in $\H^2\cup \P^1(K)$. Therefore, in the quotient by the congruence group 
$\Gamma$, the geodesic lie in Hirzebruch-Zagier divisor on the Hilbert modular surface. 
Thus, we integrate a closed algebraic differential $2$-form, (that is, a global differential $2$-form with algebraic coefficients), on the Hilbert modular surface, with boundaries Hirzebruch-Zagier divisors.

\begin{conjecture}
\label{conj period}
Let $f\in S_{k,k}(\Gamma)$ be a normalized cusp Hecke eigenform of weight $(k,k)$. Then 
\[\int_{\{p_1,p_2,p_3\}}fdz_1\wedge dz_2\]
for  geodesic triangles and
\[\int_{\{p_1,p_2;p_3,p_4\}}fdz_1\wedge dz_2\]
for geodesic diangles
are periods.
\end{conjecture}
Theorem  \ref{thm period} is a proof of Conjecture \ref{conj period} for the case of cusp form of weight $(2,2)$.

\subsection{Iteration - revisited}

We have defined iterated integrals on diangles in Definitions \ref{def it int on U type b}, \ref{def it int on U type c}. However, these definitions have to be extended to other domains of integration in order to consider iterated integrals on geodesic triangles. 

A consequence of the results from this Subsection is the following:
\begin{theorem}
\label{thm periods of it int of type c}
Iterated integrals of type {\bf{c}} on a geodesic diangle and on a geodesic triangle of algebraic differential $2$-forms on a Hilbert modular surface are periods in the sense of Kontsevich-Zagier.
\end{theorem}

Before giving the proof, we need definitions of several objects as well as their properties. In the process, we will be able to extend the definition on iterated integral on membrane when the domain of integration is a geodesic triangle.

For type {\bf{b}}, in Definition \ref{def it int on U type b}, we have a map $g:U\rightarrow \H^2$ that sends the two $\R$-foliations on $U$ into two coordinate-wise $\C$-foliations of $\H^2$. The same definition does not work when the domain $U$ is a geodesic triangle. The reason is that a geodesic triangle is either a holomorphic curve or (an anti-holomorphic curve). In both cases, a pull-back of one leaf to the geodesic triangle is a point not a line, (which is the case for the diangles).

In order to extend Definitions \ref{def it int on U type b}, \ref{def it int on U type c} to the case when the domain $U$ is a geodesic triangle, we are going to construct a new spaces using fiber product multiple times.

Now, we are going to define a space $Y_n$ associated to an iterated integral on $n$ $2$-forms on $\H^2$.
We are going to use fiber products (see \cite{Hartshorne}). 
Let $p_1$ and $p_2$ be the projections of $\H^2$ on the first and the second component, respectively.
Define $X_{ij}=\H^2$ for $1\leq i\leq n$ and $1\leq j\leq n$. (One should think of the component $X_{ij}$ as the complexification of the real coordinated $(s_i,t_j)$.) Let $C_i=\H$ for  $1\leq i\leq n$ and $C'_j=\H$ for $1\leq j\leq n$. Let 
\[X_j=X_{1j}\times_{C'_j}X_{2j}\times_{C'_j}\cdots\times_{C'_j} X_{nj}.\]
($X_j$ corresponds to the variable $t_j$)
Then
\[X_j\subset 
X_{1j}\times X_{2j}\times\cdots\times X_{nj}.\]
Let also 
\[P_j=(p_1,\cdots,p_1):X_{1j}\times X_{2j}\times\cdots\times X_{nj}\rightarrow C_1\times\cdots \times C_n\]
Let $P^\circ_j=P_j|_{X_j}$ be the restriction of $P_j$ to the subset $X_j$.
We define $Y_n$ as the fiber product of $X_1,\dots,X_n$ with respect to the morphisms $P^\circ_1,\dots,P^\circ_n$ over the base $C_1\times\cdots \times C_n$, namely
\begin{equation}
\label{def Y_n 1}Y_n=X_1\times_C\cdots\times_C X_n,
\end{equation}
where $C=C_1\times\cdots\times C_n$.
Note that $X_j$ is isomorphic to $X_{j+1}$. Let $Z_j$ be the subspace of $Y_n$ defined by setting the 
$j$- and the $(j+1)$-component of $Y_n=X_1\times_C\cdots\times_C X_n$ to be equal. (The space $Z_j$ corresponds to a boundary components obtained by letting $t_j=t_{j+1}$.)
Similarly, we could have defined $Y_n$ by defining first
\[X'_i=X_{i1}\times_{C_i}X_{i2}\times_{C_i}\cdots\times_{C_i} X_{in}\]
($X'_i$ corresponds to $s_i$)
so that 
\[X'_i\subset X_{i1}\times X_{i2}\times\cdots\times X_{in}\]
Let 
\[P'_i
=
(p_2,\dots,p_2):
X_{i1}\times X_{i2}\times\cdots\times X_{in}
\rightarrow 
C'_1\times \cdots\times C'_n\]
Define
$P'^{\circ}_i=P'_i|_{X'_i}$ to be the restriction of $P'_i$ to $X'_i$. We define $Y_n$ as the fiber product of 
$X'_1,\dots,X'_n$ with respect to the morphisms $P'^\circ_1,\dots,P'^\circ_n$ over the base $C'_1\times\cdots \times C'_n$, namely
\begin{equation}
\label{def Y_n 2}
Y_n=X'_1\times_{C'}\cdots\times_{C'}X'_n,
\end{equation}
where $C'=C'_1\times\cdots\times C'_n$.
Similarly we define $Z'_i$ to be the subspace of $Y_n$ defined by setting the 
$i$- and the $(i+1)$-component of $Y_n=X'_1\times_{C'}\cdots\times_{C'} X'_n$ to be equal.
 (The space $Z'_i$ corresponds to a boundary components obtained by letting $s_i=s_{i+1}$.)

We have given two Definitions \ref{def Y_n 1} and \ref{def Y_n 2} of the space $Y_n$. 
In the two definitions we have only exchanged the role of $p_1$ and $p_2$. We will prove that both definitions lead to the same object in the case of $n=2$. The general case is left to the reader.

\begin{lemma}
For $n=2$, the two Definitions \ref{def Y_n 1} and \ref{def Y_n 2} define isomorphic objects $Y_2$.
\end{lemma}
\proof
The space $Y_2$ can be defined as a finite limit (in a categorial sense) of a diagram in the following way.
Consider the commutative diagram

%$$
%\xymatrix{
%   &H^q(\overline{Y}_{13},F_V))
%    \ar[r] \ar[dr]
%             &H^q(\overline{Y}_{12},F_V))  \ar[dr]\\
%E_1^{*.q}: &H^q(\overline{Y}_{12,34},F_V)  \ar[ur]
%       |!{[u];[r]}\hole \ar[dr] |!{[d];[r]}\hole
%             & H^q(\overline{Y}_{23},F_V) \ar[r]
%                    &  H^q(\overline{Y}_B,F_V)\\
%  &H^q(\overline{Y}_{24},F_V)   \ar[r] \ar[ur]
%             & H^q(\overline{Y}_{34},F_V) \ar[ur]
%}
%$$

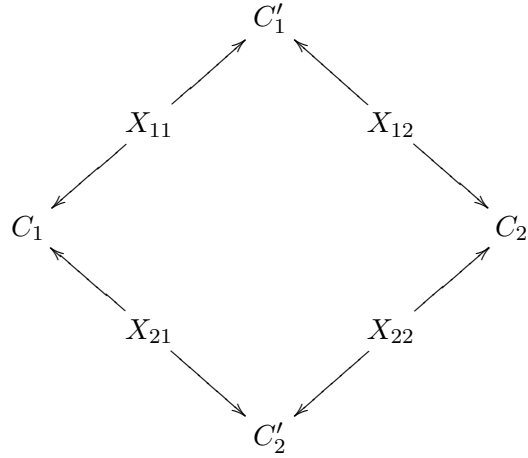
\[
\xymatrix{
         &&C'_1&&\\
         &X_{11} \ar[ur] \ar[dl] &        & X_{12}\ar[dr] \ar[ul]& \\
    C_1  &                               &&                               &C_2\\
        &X_{21} \ar[dr] \ar[ul] &        & X_{22}\ar[dl] \ar[ur]& \\
        &                      &C'_2 &&\\ }
\]

For any space $W$ such that

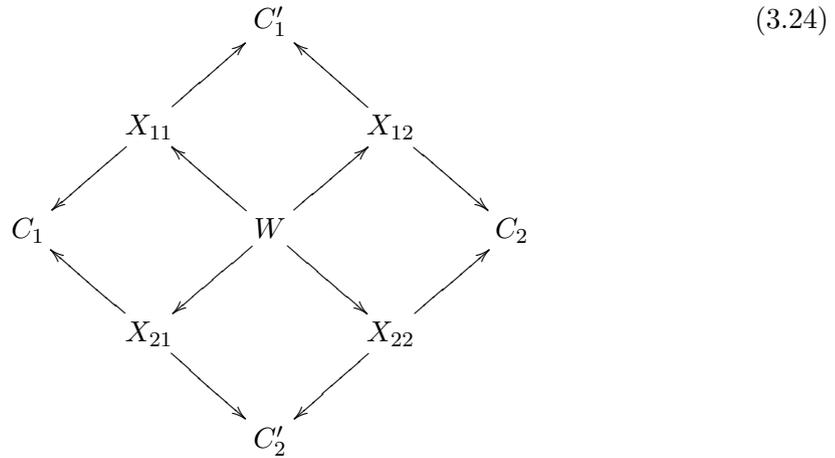
\begin{equation}
\label{diagram}
\xymatrix{
&&C'_1&&\\
         &X_{11} \ar[ur] \ar[dl] &        & X_{12}\ar[dr] \ar[ul]& \\
C_1  &                               &W  \ar[ul] \ar[ur] \ar[dr] \ar[dl]&                               &C_2\\
   &X_{21} \ar[dr] \ar[ul] &        & X_{22}\ar[dl] \ar[ur]& \\
   &                      &C'_2 &&\\ }
\end{equation}
commutes,
we have that the maps $f_{ij}:W\rightarrow X_{ij}$ factor through $g_{ij}:Y_2\rightarrow X_{ij}$ so that
$f_{ij}=g_{ij}\circ h$, for some $h:W\rightarrow Y_2$ and $Y_2$ is part of the commutative diagram

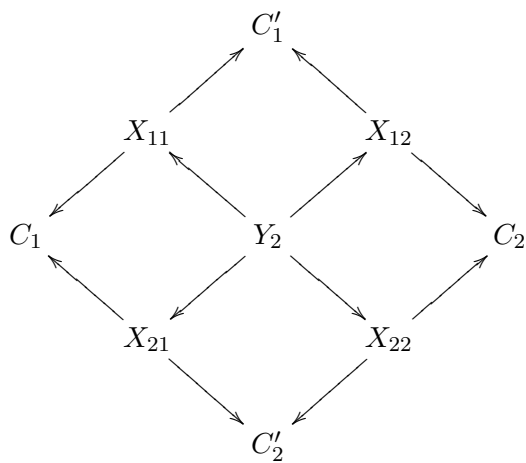
\[
\xymatrix{
&&C'_1&&\\
         &X_{11} \ar[ur] \ar[dl] &        & X_{12}\ar[dr] \ar[ul]& \\
C_1  &                               &Y_2  \ar[ul] \ar[ur] \ar[dr] \ar[dl]&                               &C_2\\
   &X_{21} \ar[dr] \ar[ul] &        & X_{22}\ar[dl] \ar[ur]& \\
   &                      &C'_2 &&\\ }
\]
In order to prove this universal property of $Y_2$ we follow the first definition of $Y_2$. It leads to the commutative diagram
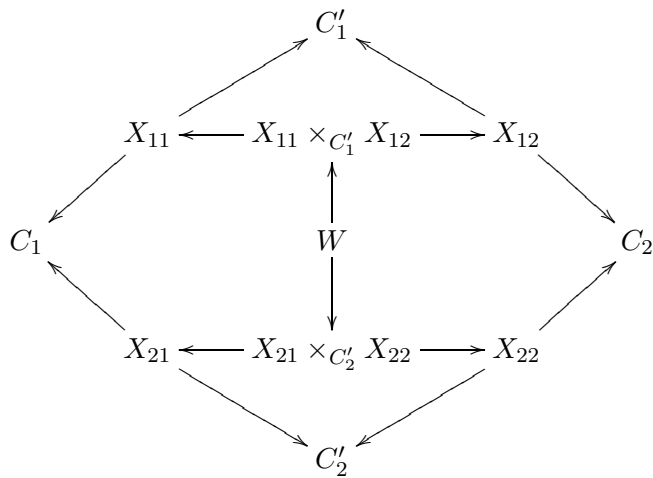
\[
\xymatrix{
&&C'_1&&\\
         &X_{11} \ar[ur] \ar[dl]&   X_{11}\times_{C'_1}X_{12}   \ar[r] \ar[l]  & X_{12}\ar[dr] \ar[ul]& \\
C_1  &                               &W \ar[u] \ar[d]&                               &C_2\\
   &X_{21} \ar[dr] \ar[ul]  &    X_{21}\times_{C'_2}X_{22}  \ar[r] \ar[l]    & X_{22}\ar[dl] \ar[ur]& \\
   &                      &C'_2 &&\\ }
\]
Then we have that  $X_1=X_{11}\times_{C'_1}X_{12}$ maps to $C=C_1\times C_2$ and also
$X_2=X_{21}\times_{C'_2}X_{22}$ maps to $C=C_1\times C_2$.
Thus the maps from $W$ to any element of the diagram factors through $Y_2=X_1\times_C X_2$.
Similarly, $W$ factors through $X'_1\times_{C'} X'_2$, where $X'_1=X_{11}\times_{C_1}X_{21}$,
$X'_2=X_{12}\times_{C_2}X_{22}$ and $C'=C'_1\times C'_2$. Since both $X_1\times_C X_2$ and 
$X'_1\times_{C'} X'_2$ are universal objects with respect to the diagram \eqref{diagram}, we have that they are isomorphic. \qed

Now, we return to the initial question of this subsection, namely, how to iterated over a geodesic triangle so that it is consistent with the current definition of iteration over a diangle. 

For an $n$-fold iteration of $n$ $2$-forms of types {\bf{b}} or {\bf{c}}, we have to specify a domain $U\subset \H^2$, $dim_\R U=2$ and a pair of permutations $\rho_1$ and $\rho_2$ of $n$ elements. We make an essential assumption that the boundary of $U\subset \H^2$, denoted by $\partial U$, projected onto the Hilbert modular surface $Y_\Gamma$ lies on a finite union of Hirzebruch-Zagier divisors. We will denote the finite union of such Hirzebruch-Zagier divisors by $HZ$.

Let 
\[P_{\rho_1,\rho_2}: 
X_{11}\times\cdots \times  X_{nn}
\rightarrow 
X_{\rho_1(1)\rho_2(1)}\times\cdots\times X_{\rho_1(n)\rho_2(n)}\]
be a projection to $n$ of the factors. Let $U_{ij}\cong U$ for $1\leq i\leq n$ and  $1\leq j\leq n$.
Let 
\[I:
U_{\rho_1(1)\rho_2(1)}\times\cdots\times U_{\rho_1(n)\rho_2(n)}
\rightarrow
X_{\rho_1(1)\rho_2(1)}\times\cdots\times X_{\rho_1(n)\rho_2(n)}\]
be the induced from the product of inclusion of the domains $U\rightarrow X$. 
We will use the following notation 
\[\underline{U}^\rho=U_{\rho_1(1)\rho_2(1)}\times\cdots\times U_{\rho_1(n)\rho_2(n)}\]
and
\[\underline{X}^\rho=X_{\rho_1(1)\rho_2(1)}\times\cdots\times X_{\rho_1(n)\rho_2(n)}.\]
Then the map $I$ becomes \[I:\underline{U}^\rho\rightarrow \underline{X}^\rho\]
Let \[J:Y_n\rightarrow \underline{X}^\rho\]
be the composition of the natural inclusion $Y_n\rightarrow X_{11}\times\cdots \times  X_{nn}$ and the projection $P_{\rho_1,\rho_2}$.
Then we define the domain of integration to be
\[\underline{U}^\rho_{Y_n}=\underline{U}^\rho\times_{\underline{X}^\rho}Y_n,\]
which is the fiber product of the maps $I$ and $J$.
Since  $I:\underline{U}^\rho\rightarrow \underline{X}^\rho$ is an inclusion, we have that the induced map
\[\underline{U}^\rho_{Y_n}\rightarrow Y_n\]  is an inclusion.

In the above constructions, we have used a parallel between type {\bf{b}} and type {\bf{c}} of iterated integrals on membranes. THe following definition allows us to extend in some sense the two types when the domain of integration is an ideal triangle
\begin{definition}
\label{def it int on U} 
(iterated integrals on membranes of types {\bf{b'}} or {\bf{c'}})
For any manifold  with corners of dimension $2$ on the Hilbert modular variety,
we define an iterated integral
\begin{equation}
\label{eq U it int}
\int_U^{\Sigma_n(\rho_1,\rho_2)}(f_1dz_1\wedge dz_2)\cdots(f_ndz_1\wedge dz_2)
=
\int_{\underline{U}^\rho_{Y_n}}J^*(f_1dz_1\wedge dz_2,\dots, f_n dz_1\wedge dz_2),
\end{equation}
where $f_kdz_1\wedge dz_2$ is a form defined on $X_{\rho_1(k)\rho_2(k)}$, for $1\leq k\leq n$.
If $Y_n$ and $\underline{U}^\rho_{Y_n}$ are constructed in the setting of type {\bf{b}} iterated integrals on membranes, then the above definition is of iterated integral on membranes of type {\bf{b'}}. 
Similarly, if $Y_n$ and $\underline{U}^\rho_{Y_n}$ are constructed in the setting of type {\bf{c}} iterated integrals on membranes, then the above definition is of iterated integral on membranes of type {\bf{c'}}. 
\end{definition}

If $U$ is a diangle, then the relation of the above integral to the ones defined by iterated integrals over membranes is the following. The integral 
\[\int_{U}^{\Sigma_n(\rho_1,\rho_2)}(f_1dz_1\wedge d z_2)\cdots(f_ndz_1\wedge dz_2)\] 
is the sum of the integrals from Definitions \ref{def it int on U type b} or \ref{def it int on U type c}, namely, the sum 
\[\sum_{\rho\in \Sigma_n}\int_{U}^{(\rho\rho_1,\rho\rho_2)}(f_1dz_1\wedge d z_2)\cdots(f_ndz_1\wedge dz_2)\]
over the orbit of the diagonal action of the permutation group $\Sigma_n$ 
on any chosen pair of  permutations $(\rho_1,\rho_2)$.

\begin{proposition} (Properties of the iterated integral \ref{eq U it int})
(a) the iterated integral \ref{eq U it int} is well-defined when $U$ is an ideal triangle both for types {\bf{b}} and {\bf{c}};

(b) the iterated integral \ref{eq U it int} for type {\bf{c}} is a period if $U$ is an ideal triangle or an ideal diangle, when $f_1,\dots,f_n$ are normalized Hecke eigenforms of weight $(2,2)$.

(c) the iterated integral \ref{eq U it int}, both for types {\bf{b}} and {\bf{c}},
 is homotopy invariant with respect to a homotopy that varies within the  divisors \[J^{-1}(X_{\rho_1(1)\rho_2(1)}\times\cdots\times p_1^{-1}(q_i)\times\cdots \times X_{\rho_1(n)\rho_2(n)}),\]
where $q_i$ is a point of  $X_{\rho_1(i)\rho_2(i)}$ for fixed $i$ and 
$p_1:X_{\rho_1(i)\rho_2(i)}\rightarrow C$;
or a homotopy that varies within the  divisors 
\[J^{-1}(X_{\rho_1(1)\rho_2(1)}\times\cdots\times p_2^{-1}(q_i)\times\cdots \times X_{\rho_1(n)\rho_2(n)}),\]
where $q_i$ is a point of  $X_{\rho_1(i)\rho_2(i)}$ for fixed $i$ and 
$p_2:X_{\rho_1(i)\rho_2(i)}\rightarrow C'$.
\end{proposition}
\proof (a) The integral \ref{eq U it int} is well defined for any two dimensional submanifold with corners of the Hilbert modular variety. (\cite{BS})

(b) The iterated integral \ref{eq U it int} is a period since:

(1)  a Hilbert modular variety can be defined over a number field 

(2) the normalized Hecke eigenforms $f_1,\dots,f_n$ of weight (2,2) can be realized as algebraic differential forms on the Hilbert modular variety;

(3) the boundary of the region of integration $\overline{U}^\rho_{Y_n}$ is a divisor on $Y_n$, namely,
\[\bigcup_{i=1}^{n}HZ_i,\]
where 
\[HZ_i= J^{-1}(X_{\rho_1(1)\rho_2(1)}\times\cdots\times HZ\times\cdots \times X_{\rho_1(n)\rho_2(n)})\]
is a divisor of $Y_n$ obtained by a pull-back of a divisor whose $i$-th component is a Hirzebruch-Zagier
divisor $HZ$ and the rest of the factors are $X_{\rho_1(k)\rho_2(k)}$ for $k\neq i$.

(c) The proof is essentially the same as the one of Theorem \ref{thm h invariance}.
\qed

The domain $\underline{U}_{Y_n}$ might be cut into disconnected components by the divisors $Z_i$'s and $Z'_j$'s. In order to choose a connected component we need to define another region of integration. 
Recall: For the case of intreated integrals on membranes of type {\bf{b}},
 $p_1:\H^2\rightarrow C$ and $p_2:\H^2\rightarrow C'$ be projections onto the first and the second component with $C\cong \H$ and $C'\cong \H$.
 
For the case of iterated integrals on membranes of type {\bf{c}},
$p_1=\alpha_1\circ\pi$ and $p_2=\alpha_2\circ\pi$ 
are compositions of the map from the universal cover to the Hilbert modular surface 
\[\pi:\H^2\rightarrow X_\Gamma\] and \[\alpha_1,\alpha_2:X_\Gamma\rightarrow \P^1\] 
be two algebraically independent rational functions on the Hilbert modular surface
and $C_i\cong\P^1$ and $C'_j\cong \P^1$ for $1\leq i\leq n$ and $1\leq j\leq n$.

Let $q_0,q_1,r_0,r_1\in \P^1$ be points. Let each of $Q_0,Q_1,R_0$ and $R_1$ be a connected component of $p_1^{-1}(q_0)$, of $p_1^{-1}(q_1)$, of $p_2^{-1}(r_0)$, and of $p_2^{-1}(r_1)$,
respectively.

Let $V\rightarrow \H^2$ be  a domain in $\H^2$ such that its boundaries lie on the union
\[
Q_0\cup Q_1\cup R_0\cup R_1,
\]
and its interior does not meet the union $Q_0\cup Q_1\cup R_0\cup R_1.$
We define the divisors $Z_0,Z_n,Z'_0,Z'_n$  of $Y_n$, which will have the following meaning: 
$Z_0$ will be the beginning of the integration of the $t_1$ variable ($t_1=0$), 
$Z_n$ will be the end of the integration of the $t_n$ variable  ($t_n=1$),
$Z'_0$ will be the beginning of the integration of the $s_1$ variable ($s_1=0$),
and
$Z'_n$ will be the end of the integration of the $s_n$ variable  ($s_n=1$).
We define them as the a fiber product 
\[Z_0=Q_0\times_C X_2\times_C \cdots\times_C X_n,\]
\[Z_n=X_1\times_C\cdots\times_C X_{n-1}\times_C Q_1,\]
\[Z'_0=R_0\times_{C'} X'_2\times_{C'}  \cdots\times_{C'}  X'_n\]
and
\[Z'_n=X'_1\times_{C'} \cdots\times_{C'}  X'_{n-1}\times_{C'}  R_1.\]

We will use the following notation 
\[\underline{V}^\rho=V_{\rho_1(1)\rho_2(1)}\times\cdots\times V_{\rho_1(n)\rho_2(n)}\]
and
\[\underline{X}^\rho=X_{\rho_1(1)\rho_2(1)}\times\cdots\times X_{\rho_1(n)\rho_2(n)}.\]
Then the map $I'$ becomes \[I':\underline{V}^\rho\rightarrow \underline{X}^\rho\]
Let \[J':Y_n\rightarrow \underline{X}^\rho\]
be the composition of the natural inclusion $Y_n\rightarrow X_{11}\times\cdots \times  X_{nn}$ and the projection $P_{\rho_1,\rho_2}$.
Then we define the domain of integration to be
\[\underline{V}^\rho_{Y_n}=\underline{V}^\rho\times_{\underline{X}^\rho}Y_n,\]
which is the fiber product of the maps $I'$ and $J'$.
Since  $I':\underline{U}^\rho\rightarrow \underline{X}^\rho$ is an inclusion, we have that the induced map
\[\underline{V}^\rho_{Y_n}\rightarrow Y_n\]  is an inclusion.

Then the divisors $Z_0,Z_1,\dots,Z_{n-1},Z_n$ and $Z'_0,Z'_1,\dots,Z'_{n-1},Z'_n$ cut out from 
$\underline{V}_{Y_n}$ a product of two  $n$-simplices, corresponding to the region where $\{0\leq s_1\leq \cdots\leq s_n\leq 1\}\times \{0\leq t_1\leq \cdots\leq t_n\leq 1\}$ is embedded. 
Denote by $\overline{V}^\rho_{Y_n}$ the connected components of 
$\overline{V}^\rho_{Y_n}$
that contains the image of
$\{0\leq s_1\leq \cdots\leq s_n\leq 1\}\times \{0\leq t_1\leq \cdots\leq t_n\leq 1\}$
under the map $g$ from Definition \ref{def membranes 1c}.

\vspace{.3cm}
\proof (of Theorem \ref{thm periods of it int of type c})
We consider the type of iterated integrals defined in  Definition \ref{def it int on U type c}.
Using the above notation the domain of integration is $U$, where $U\subset V$. We define
\[\overline{U}^\rho_{Y_n}=\underline{U}^\rho_{Y_n}\cap \overline{V}^\rho_{Y_n}\]
Then the boundary of  $\overline{U}^\rho_{Y_n}$ lies on the union of divisors
\[\partial \overline{U}^\rho_{Y_n}\subset \left(\bigcup_{i=1}^nZ_i\right)\cup\left(\bigcup_{j=1}^nZ'_j\right).\]
The normalized Hecke eigenforms of weight (2,2) can be realized as algebraic differential forms on the Hilbert modular variety. Then the iterated integrals on a membrane  of type {\bf{c}} over the domain $U$ is a period since: 

(1)  a Hilbert modular variety can be defined over a number field 

(2) the normalized Hecke eigenforms $f_1,\dots,f_n$ of weight (2,2) can be realized as algebraic differential forms on the Hilbert modular variety;

(3) the boundary of the region of integration $\overline{U}^\rho_{Y_n}$ is a divisor on $Y_n$, namely,
\[\left(\bigcup_{i=1}^nZ_i\right)\cup\left(\bigcup_{j=1}^nZ'_j\right).\]
\qed

\subsection{Generating series and relations}

In this Subsection we examine the generating series of iterated integrals on membranes (of types {\bf{b'}} or {\bf{c'}}), evaluated at geodesic triangles and geodesic diangles. We prove relations among them. Most importantly, the generating series $J$ will be used in Subsection 3.5 for Defining non-commutative Hilbert modular symbols. Moreover, the relations that we prove in this Section, will be interpreted as cocycles or as coboundaries of the the non-commutative Hilbert modular symbols satisfy in Subsection 3.5. 

\begin{definition}
Let $f_1,\dots,f_m$ be $m$ cusp forms with respect to a Hilbert modular group $\Gamma$.
Let $f_1dz_1\wedge dz_2,\dots,f_mdz_1\wedge dz_2$ be the corresponding differential forms, defining the generating series.
Let $J(p_1,p_2,p_3)$ be the generating series $J$ evaluated at the geodesic triangle with vertices $p_1,p_2,p_3$.
Let $J(p_1,p_2;p_3,p_4)$ be the generating series $J$ evaluated at the geodesic diangle $\{p_1,p_2;p_3,p_4\}$. 
\end{definition}
Both $J(p_1,p_2,p_3)$ and $J(p_1,p_2;p_3,p_4)$ will be called non-commutative Hilbert modular symbols after the action of the arithmetic group is included (see Definition \ref{def c1 c2}).

\begin{theorem}
\label{thm NC relations}
The generating series $J$ is one of the following types {\bf{b}}, {\bf{c}}, {\bf{b'}} or {\bf{c'}}. Note that
$J(p_1,p_2;p_3,p_4)$ is defined for all the types, while $J(p_1,p_2,p_3)$ is defined only for types {\bf{b'}} or {\bf{c'}}. Then the generating series $J(p_1,p_2,p_3)$ and $J(p_1,p_2;p_3,p_4)$ satisfy the following relations:

1. If $\sigma$ is a permutation of the set $\{1,2,3\}$ then  
\[J(p_{\sigma(1)},p_{\sigma(2)},p_{\sigma(3)})
=
J^{sign(\sigma)}(p_1,p_2,p_3).\] 

2. If $p_1,p_2,p_3,p_4$ are four points on the same holomorphic (or anti-holomorphic) curve of the type $\gamma^*\Delta$ then
\begin{align*}
1=&J(p_1,p_2,p_3)J(p_2,p_3,p_4)\\
&J(p_2,p_1,p_4)  J(p_1,p_4,p_3)\\
\end{align*}
and

 3. If $p_1,p_2,p_3,p_4$ are four points on the same holomorphic (or anti-holomorphic) curve of the type $\gamma^*\Delta$ then
\[1=J(p_1,p_2;p_3,p_4).\]

4. For every four points $p_1,p_2,p_3,p_4$, we have the following relation based on the orientation of the domain
 \begin{align*}
J(p_2,p_1;p_3,p_4)=&J(p_1,p_2;p_4,p_3)=\\
=&J^{-1}(p_2,p_1;p_4,p_3)=\\
=&J^{-1}(p_1,p_2;p_3,p_4).
\end{align*}

5. For every five points $p_1,p_2,p_3,p_4,p_5$, we have 
 \[J(p_1,p_2;p_3,p_4) J(p_1,p_2;p_4,p_5) = J(p_1,p_2;p_3,p_5).\]

6.  For every four points $p_1,p_2,p_3,p_4$, we have the following relation, based on the boundary of a thickened tetrahedron,
\begin{align*}
1=&J(p_1,p_2,p_3)J(p_2,p_3,p_4)\\
&J(p_2,p_1,p_4)  J(p_1,p_4,p_3)\\
&J(p_1,p_2;p_3,p_4) J(p_2,p_3;p_1,p_4) J(p_3,p_1;p_2,p_4) \\
&J(p_3,p_4;p_1,p_2)  J(p_1,p_4;p_2,p_3) J(p_2,p_4;p_3,p_1).
\end{align*}
\end{theorem}
\proof
For part 1, let $\sigma$ be an odd permutation. Let $U$ be an union of two triangles along one of their the edges. 
Let the first triangle be with vertices $p_1,p_2,p_3$ and the second triangle be with vertices $p_3,p_2,p_1$ with the opposite orientation. We can glue the two triangles along the edge $p_1p_2$. (Glueing along any other edge would lead to the same result for the corresponding generating series.)
From the shuffle product formula Theorem \ref{thm shuffle} (iii), 
it follows that $J(U)=J(p_1,p_2,p_3)J(p_3,p_2,p_1)$.  (Note that the product is not the product in the ring $R$. It is induced by a shuffle product of iterated integrals on membranes.)
From the second homotopy invariance theorem 
(Theorem \ref{thm h-invariance II}) 
it follows that the generating series $J(U)$ depends on $U$ up to homotopy, 
which keeps the boundary components $p_2p_3$, $p_3p_2$, $ p_1p_3$  and $p_3p_1$ on fixed union holomorphic curves. 
We can contract $U$ to its boundaries $\partial U$
so that the contracting homotopy ``keeps the boundary components" 
on a fixed union of holomorphic curves. 
Therefore, $J(U)=J(\partial U)=1$. 

Parts 2, 4 and 5 can be proven similarly. 

For part 3, if $p_1,p_2,p_3,p_4$ belong to the same holomorphic (or anti-holomorphic) curve then the corresponding diangle has no interior, since the two edges will coincide. Recall that the edges of the diangle are defined via unique geodesic triangles lying on a holomorphic (or anti-holomorphic) curve. 

The proof of part 6 is essentially the same as the one for part 1; however, we will prove it independently, since it is a key property of the non-commutative Hilbert modular symbol. Consider a thickened tetrahedron with vertices $p_1,p_2,p_3,p_4$. The faces of the thickened tetrahedron are precisely the ones listed in the product of property 6. The whole product is equal to $J(V)$, where $V=\mbox{\it{union of all faces of the thickened tetrahedron}}$. 
From the second homotopy invariance theorem it follows that the generating series $J(V)$ depends on $V$ up to homotopy, which keeps the boundary components on a fixed union holomorphic curves. 
Since $V$ bounds a contractible $3$-dimensional region (a thickened tetrahedron), from Theorem
\ref{thm h-invariance II}, it follows that $J(V)=J(point)=1$. \qed

\subsection{Definition of non-commutative Hilbert modular symbols}
In this Subsection we define non-commutative Hilbert modular symbols. They are analogues of Manin's non-commutative modular symbol (see \cite{Man}), applicable to the Hilbert modular group. Instead of iterated path integrals that Manin uses, we use a higher dimensional analogue, defined in Section 2.

Usually, a modular symbol represents a cohomology class. Manin's non-commutative modular symbol represent a non-commutative $1$-st cohomology class. We would like to say that the non-commutative Hilbert modular symbols represent non-commutative cohomology classes, which we formulate as Conjectures.

After defining the non-commutative Hilbert modular symbols, we prove some of their properties. These properties will be interpreted  intuitively as co-cycle conditions or as co-boundary conditions. The approach in this Subsection is more geometric. The purpose for presenting them here is to give many examples of relations and to help establishing a suitable cohomology theory that would truly capture these relations in a more structured way.

The cocycle interpretation is only for intuition it is not precise. The formula holds for geometric reasons. Note that the composition is not the multiplication in the ring R, it is given by the shuffle product (see Theorem \ref{thm shuffle}), which works for the generating series on iterated integrals on membranes.  The multiplication is written linearly as we would multiply several elements in a group or in a ring; however, the multiplication is two-dimensional among regions with common boundaries. 

In the next Subsection will give some intuition about higher categories for the purpose of giving more structure to the non-commutative Hilbert modular symbols and for a possible approach to defining such a $1$-st and $2$-nd non-commutative cohomology class.

For definitions of iterated integrals on membranes of types  {\bf{b}}, {\bf{c}}, {\bf{b'}} and {\bf{c'}} see Definitions \ref{def membranes 1b}, \ref{def membranes 1c}) for types  {\bf{b}}, {\bf{c}}, and 
\ref{def it int on U} for types  {\bf{b'}} and {\bf{c'}}.
\begin{definition}
\label{def c1 c2}
We define non-commutative Hilbert modular symbols as a generating series of iterated integrals on membranes of types {\bf{b}}, or {\bf{c}}, or {\bf{b'}}, or {\bf{c'}} over a geodesic diangle by
\[c^1_{p_1,p_2;p_3}(\gamma)=J(p_1,p_2;p_3,\gamma p_3).\]
We also define non-commutative Hilbert modular symbols as a generating series of iterated integrals on membranes of types {\bf{b'}} or {\bf{c'}} over a geodesic triangle by
\[c^2_{p}(\gamma,\delta)=J(p,\gamma p, \gamma \delta p),\]
where $p,p_1,p_2,p_3$ are cusp points in $\H^2\cup\P^1(K)$ and 
$\beta, \gamma, \delta \in SL_2({\cal{O}}_K)$.
\end{definition}

We are going to define an action of 
$Mat_2({\cal{O}}_K)^+$
on the generating series $J^c$, where $Mat_2({\cal{O}}_K)^+$ is the semi-group of  $2\times  2$ matrices with a totally positive determinant.

In order to interpret $c^1(\gamma)$ and $c^2(\gamma,\delta)$ as cocycles, we are going to define an action of the semi-group $Mat_2({\cal{O}}_K)^+$ on the whole ring $R$, where the generating series take values. Such an action can be given via Hecke operators. 

For simplicity, we shall assume that ${\cal{O}}_K$ has narrow class number $1$. We consider all Hecke eigenforms of weight $(2,2)$ with respect to  $Mat_2({\cal{O}}_K)^{+}$. Now, let $u$ be a unit such that $u_1>0$ and $u_2<0$, where $u_1$ and $u_2$ are the images of $u$ under the two real embeddings of $K$ into $\R$. It exists, since the narrow class group is trivial. 
(For example $K=\Q(\sqrt{2})$ is such a field.)
We define an action of $\gamma\in Mat_2({\cal{O}}_K)$ on the ring $R$ (Definition \ref{def ring R}), where the generating series takes values. We 
define \[\gamma\bullet f\mapsto T_\gamma (f)\] if $\gamma\in Mat_2({\cal{O}}_K)^{+}$.
Let $f_1,\dots, f_m$ be a basis of Hecke eigenforms of the space of the cusp form of weight $(2,2)$. To each $f_i$ we associate two of the generators of $R$ $X_i$ and $Y_i$.
Let $X_1,Y_1,\dots,X_m,Y_m$ be generators of $R$. Then the action of $\gamma\in Mat_2({\cal{O}}_K)^{+}$ is given by  $T_\gamma(X_i)=c(\gamma,f_i)X_i$
and $T_\gamma(Y_i)=Y_i$, where $c(\gamma,f_i)$ the eigenvalue of the Hecke operator.

 In this setting the group action, namely, the action of the Hilbert modular group is trivial. This trivial action extend to the the action of $T_1=id$ on the whole ring $R$. In fact, for an element $\beta\in SL_2({\cal{O}}_K)$, the trivial action on $c^1_{p_1,p_2;p_3}$
and on $c^2_{p}$ can be realized as follows:
\[(\beta c^1_{p_1,p_2;p_3})(\gamma)=c^1_{\beta p_1, \beta p_2; \beta  p_3}(\beta\gamma)\]
and 
\[(\beta c^2_{p})(\gamma,\delta)=c^2_{\beta p}(\beta p, \beta \gamma p, \beta \gamma \delta p).\]
The last two relations hold, since for a cusp form of weight $(2,2)$ the differential form $fdz_1\wedge dz_2$ is invariant under the action of the Hilbert modular group $\Gamma$. Algebraically, for any geodesic diangle, we have
\[\beta J(p_1,p_2;p_3,p_4)=J(p_1,p_2;p_3,p_4)=J(\beta p_1,\beta p_2;\beta p_3,\beta p_4).\]
similarly for a geodesic triangle,
\[\beta J(p_1,p_2,p_3)=J(p_1,p_2,p_3)=J(\beta p_1,\beta p_2,\beta p_3).\]

The relations among the symbols are based on two properties: composition via shuffle product Theorem \ref{thm shuffle} (iii) and the homotopy invariance (Theorems \ref{thm h invariance} and \ref{thm h-invariance II}).

\begin{conjecture}
\label{thm c1 c2}
The non-commutative Hilbert modular symbol $c^1_{p_1,p_2;p_3}$ is a $1$-cocycle. Moreover, if we change the point  $p_3$ to $q_3$, then the cocycle changes by a coboundary.
\end{conjecture}
Property 5 of Theorem \ref{thm NC relations} can be interpreted as a $1$-cocycle relation. 
Consider the analogy with non-commutative $1$-cocycle of a group acting on a non-commutative ring, we define the boundary of $c^1_{p_1,p_2;p_3}$ by
\[
dc^1_{p_1,p_2;p_3}(\beta,\gamma)
=c^1_{p_1,p_2;p_3}(\beta)(\beta c^1_{p_1,p_2;p_3})(\gamma)(c^1_{p_1,p_2;p_3}(\beta\gamma))^{-1}\]
The action of $\beta$ on the cocycle is given in Definition \ref{def c1 c2}.
In contrast to a $1$-st non-commutative cocycle (see for example Kenneth Brown \cite{Br}), here we have two-dimensional composition of symbols, that is, one can compose the symbols as two-morphisms in a two-category.

Then
\begin{align}
\nonumber dc^1_{p_1,p_2;p_3}(\beta,\gamma)
&=J(p_1,p_2;p_3,\beta p_3)(\beta  J(p_1,p_2;p_3,\gamma p_3))J^{-1}(p_1,p_2;p_3,\beta\gamma p_3)\\
\nonumber&=
J(p_1,p_2;p_3,\beta p_3)J(p_1,p_2;\beta p_3,\beta\gamma p_3)J^{-1}(p_1,p_2;p_3,\beta\gamma p_3)\\
&=
1.
\end{align}

If we change $p_3$ to $q_3$ then the cocycle changes by a coboundary.
Let $b^0=J(p_1,p_2;p_3,q_3)$ be a $0$-cochain. Then 
\begin{align}
\nonumber c^1_{p_1,p_2;q_3}(\gamma)
&=
J(p_1,p_2;p_3,\gamma p_3)\\
\nonumber &=J(p_1,p_2;p_3,q_3)J(p_1,p_2;q_3,\gamma q_3)J(p_1,p_2;\gamma q_3,\gamma p_3)\\
\nonumber &=J(p_1,p_2;p_3,q_3)J(p_1,p_2;q_3,\gamma q_3)(\gamma  J(p_1,p_2;p_3,q_3))^{-1}\\
&=b^0 c^1_{p_1,p_2;q_3}(\gamma)(\gamma b^0)^{-1}
\end{align}

\begin{conjecture}
\label{thm c2}
The non-commutative Hilbert modular symbol $c^2_{p}(\beta,\gamma)$ satisfies a $2$-cocycle relation. Moreover, if we change the point $p$  to $q$, then the cocycle changes by a coboundary up to terms involving $c^1$.
\end{conjecture}
Recall \[c^2_{p}(\beta,\gamma)=J(p,\beta p, \beta \gamma p).\]
Then $c^2_p$ satisfies a  $2$-cocycle condition up to a multiple of the $1$-cocycle $c^1_{q_1,q_2;q_3}$ for various points $q_1,q_2,q_3$. For the $2$-cocycle relation, we compute $dc^2_p(\beta,\gamma,\delta)$.
\begin{align}
\label{eq dc2}
\nonumber dc^2_p(\beta,\gamma,\delta)=&c^2_p(\beta,\gamma)c^2(\beta,\gamma\delta)(c^2(\beta\gamma,\delta))^{-1}(\beta\cdot c^2(\gamma,\delta))^{-1}=\\
\nonumber
=&J(p,\beta p,\beta\gamma p)J(p,\beta p,\beta\gamma\delta p)\times\\
&\times J(p,\beta\gamma p, \beta\gamma\delta p)^{-1}J(\beta p, \beta\gamma p, \beta\gamma\delta p)^{-1}
\end{align}

In order to have $dc^2_p(\beta,\gamma,\delta)=1$, we must multiply by suitable values of $c^1$, corresponding to edges of a certain thickened tetrahedron. Then
\begin{align*}
dc^2_p(\beta,\gamma,\delta)\times&
[
c^1_{p,\beta p; \beta\gamma p}((\beta\gamma)\delta(\beta\gamma)^{-1})c^1_{\beta p, \beta\gamma p;p}(\beta\gamma\delta)
c^1_{\beta\gamma p,p;\beta p}((\beta)\gamma\delta\beta^{-1})\times\\
&\times c^1_{\beta \gamma p, \beta\gamma\delta p;p}(\beta)
c^1_{p, \beta\gamma\delta p;\beta p}(\beta\gamma\beta^{-1})
c^1_{\beta p,\beta\gamma\delta p;\beta\gamma p}((\beta\gamma)^{-1})
]=\\
=&[
c^2_p(\beta,\gamma)c^2(\beta,\gamma\delta)(c^2(\beta\gamma,\delta))^{-1}(\beta\cdot c^2(\gamma,\delta))^{-1}
]
\times\\
&\times
[
c^1_{p,\beta p; \beta\gamma p}((\beta\gamma)\delta(\beta\gamma)^{-1})c^1_{\beta p, \beta\gamma p;p}(\beta\gamma\delta)
c^1_{\beta\gamma p,p;\beta p}((\beta)\gamma\delta\beta^{-1})\times\\
&\times c^1_{\beta \gamma p, \beta\gamma\delta p;p}(\beta)
c^1_{p, \beta\gamma\delta p;\beta p}(\beta\gamma\beta^{-1})
c^1_{\beta p,\beta\gamma\delta p;\beta\gamma p}((\beta\gamma)^{-1})
]=\\
=&[J(p,\beta p,\beta\gamma p)J(p,\beta p,\beta\gamma\delta p)\times\\
&\times J(p,\beta\gamma p, \beta\gamma\delta p)^{-1}J(\beta p, \beta\gamma p, \beta\gamma\delta p)^{-1}]\times\\
&\times[
J(p,\beta p;\beta \gamma p,\beta\gamma\delta p) 
J(\beta p,\beta\gamma p;p,\beta\gamma\delta p) 
J(\beta\gamma p,p;\beta p,\beta \gamma\delta p)\times \\
&\times J(\beta\gamma p,\beta\gamma \delta p;p,\beta p)  
J(p,\beta\gamma\delta p;\beta p,\beta\gamma p) 
J(\beta p, \beta\gamma \delta p; \beta\gamma p,p)
]=\\
=&1.
\end{align*}
The first equality follows from Equation \eqref{eq dc2}. The second equality follows from the definition of the symbols.  And the last equality follows from Property 6 of Theorem \ref{thm NC relations} 
with $(p_1,p_2,p_3,p_4)=(p,\beta p, \beta\gamma p, \beta,\gamma\delta p)$.  Therefore, we obtain that
$dc^2_p(\beta,\gamma,\delta)$ is $1$
up to values of the $1$-cocycle $c^1$.

Now we would like to have that $c^2_p$ and $c^2_q$ are homologous. 
\begin{conjecture} The conjectural cocycles $c^2_p$ and $c^2_q$ are homologous 
\[c^2_p(\beta,\gamma)=c^2_q(\beta,\gamma)[db_{pq}^1(\beta,\gamma)]\prod_i J(D_i),\]
up to a product of $J(D_i)$,where $D_i$ are geodesic diangles.
\end{conjecture}

Before we proceed, we would like to make an analogy between 1-dimensional and 2-dimensional cocycles. For the 1-dimensional cocycle, the property that it is a cocycle uses the geometry of a triangle where the faces of the triangle are essentially the 1-cocycle. We want commutativity of the triangular diagram. We think of the commutativity of the diagram as follows: consider the interior of the triangle as a homotopy of paths and we think of the $1$-cocycle as a homotopy invariant function.
For the $2$-cocycle, the $2$-cocycle relation is represented by the faces of a tetrahedron. By a 'commutativity' of the diagram, we mean a homotopy invariant $2$-cocycle and a homotopy from one of the faces to the union of the other three faces.

The comparison that $c^1_{p_1,p_2;p_3}$ and $c^1_{p_1,p_2;q_3}$ are homologous is given by a square-shaped diagram. The analogy with dimension $2$ is that the cocycles $c^2_p$ and $c^2_q$ are two faces of an octahedron. The vertices associated to $c^2_p(\beta,\gamma)$ are $(p,\beta p,\beta\gamma p)$ and the vertices associated to $c^2_q$ are $(q,\beta q,\beta\gamma q)$. 
The two faces will be opposite to each other on the octahedron $Oct$ so that the three pair of opposite vertices are $(p,\beta \gamma q)$, $(\beta p, q)$ and $(\beta\gamma p ,\beta q)$. The remaining $6$ faces are combined into two triples. Each of them corresponds to a coboundary of a 1-chain.

Let \[b^1_{p,q}(\beta)=[J(p,q,\beta p)J(q,\beta q,\beta p)][J(q,\beta p;p,\beta q)].\]
Consider the action of $\gamma\in\Gamma$ on $b^1$, by action on each point in the argument of $J$, denoted as before by $\gamma\cdot b^1$. Then, we define
\[db^1_{p,q}(\beta,\gamma)
=
b^1_{p,q}(\beta)[\beta\cdot b^1_{p,q}(\gamma)][b^1_{p,q}(\beta\gamma)]^{-1},
\]
where $\beta\cdot b^1_{p,q}(\gamma)=[J(\beta p,\beta q,\beta\gamma  p)J(\beta q,\beta \gamma q,\beta \gamma  p)][J(\beta q,\beta \gamma p;p;\beta \gamma q)]$.

Consider the above octahedron $Oct$. 
Remove from it the tetrahedron $T$ with vertices $(p,q,\beta\gamma q, \beta\gamma p)$.
Then the triangles of the remaining geometric figure are precisely the triangles in the definitions of $c^2_p(\beta,\gamma)$, $c^2_q(\beta,\gamma)$ and $db^1_{p,q}(\beta,\gamma)$.
Now, consider thickening of the edges, which are common for two triangles. It can be done in the following way. Instead of any triangle, we can take a geodesic triangle. The two triangles that had a common edge might have only two common vertices. Then the region between the two geodesic, one for each of the geodesic triangles, forms the induced diangle. Take $J$ of the induces diangles from the octahedron $Oct$ and $J^{-1}$ of the induced diangles from the tetrahedron $T$.  Their product gives $\prod_i J(D_i)$. The equality holds because we apply $J$ to the union of the faces of the thickened $Oct-T$, which gives $1$.

%The Hilbert modular symbols satisfy non-commutative $1$-cocycle and $2$-cocycle conditions. Unlike a non-commutative group, where the multiplication is linear, non-commutative $2$-cocycles (in terms of $J$) multiply according to the relative locations of the $2$-dimensional  
%regions of integration via shuffle product. The relations among the symbols are based on two properties: composition via shuffle product and the homotopy invariance.
%Note that from the shuffle product from Theorem \ref{thm shuffle} (iii), we have that the product uses the location of the regions of integration. In particular, the shuffle product from Theorem \ref{thm shuffle} (iii) provides a proof of Lemma \ref{lemma product of paths}, which generalizes compositions of iterated path integrals over a concatenation of paths, which is a key ingredient in the non-commutative $1$-cocycle that Manin defines to be the non-commutative modular symbol.

\subsection{A two-category}
Why do we need a two-category? Is there an example of a sheaf on this category/topology? How does the non-commutative Hilbert modular symbols represents a sheaf? 

The ideas presented in this Subsection will be developed in a follow-up paper. Here we present the basic constructions that give justification for the conjectures that the non-commutative Hilbert modular symbols $c^1$ and $c^2$ are co-cycles in some categoric and sheaf-theoretic setting. For sheaves on $2$-categories one may consult Street \cite{St}. Since our $2$-morphisms are invertible one may also use Lurie's \cite{L} constructions of sheaves on higher categories.

We are going to construct a $2$-category $C$ and a sheaf $J$ on the $2$-category $C$. 
We define $p$ to be an object of the $2$-category $C$ if $p$ is a cusp point, that is $p\in \P^1(K)$.
We define $1$-morphisms in the following way. 
Let $\sigma$ be the geodesic connecting $0$ and $\infty$ that lies on the diagonal 
$\Delta=i(\H)\subset \H\times\H$. 
There is unique such geodesic. 
All geodesics $\gamma^*\sigma$ together with a choice of orientation are defined to be $1$-morphisms, where $\gamma\in PGL_2(K)$. 
We define the $1$-morphisms of $C$ to be a finite concatenation of geodesics of the type 
$\gamma^*\sigma$ or the trivial path whose image coincides with a cusp point.
Consider ideal triangles and ideal diangles as cells from which we build manifolds with corners. 
A $2$-morphism is a finite union of manifolds with corners, made from finitely many ideal dangles and ideal triangles, which is path connected and has orientation.
 
The boundary of a $1$-morphism is a union of two object objects - the starting point and the ending point of the directed path. The boundary of a $2$-morphism (a $2$-manifold with corners) is a finite union of $1$-morphisms (oriented loops), where the orientation of the loops on the boundary  is induced by  the orientation of the $2$-manifold with corners. 

Now we are going to define a $2$-sheaf $J$, whose values on a $2$-morphism will be in a subset of the ring $R$ and whose values of an  object and on a $1$-morphism will be a subset of a countable product of the ring $R$ with itself. 

As always, $S_{2,2}(\Gamma)$ denotes the space  of cusp forms of weight $(2,2)$ with respect to the group $\Gamma$. Here we will consider this space as the space of holomorphic $2$ forms on $\H\times\H$, which vanish on the cusps and which can dessend to the Hilbert modular surface $X_\Gamma=\Gamma\backslash(\H^2\cup \P^1(K))$.
Every $n$-tuple of such holomorphic forms $\Omega\in (S_{2,2}(\Gamma))^n$ 
defines a value on a $2$-morphism $f$ in $C$. Let this values be the generating series $J_f(\Omega)$. 
Let $J_f$ be the collection of all values $J_f(\Omega)$ for all
 $\Omega\in (S_{2,2}(\Gamma))^n$. 
Let $e$ be a $1$-morphism. We say that $e$ is in the boundary of a $2$-morphism $f$, 
denoted by $e\subset \partial f$, 
if the image of the loop $e$ is in the boundary of the image of the membrane $f$ 
together with the induced orientation on $e$ from $f$. We say that an object $p$ is in the boundary of a $1$-morphism $e$, 
denoted by $p\in \partial e$, if $p$ is a source or a target of $e$.
We define the values of $J$ on a $1$-morphism $e$ to be
the product 
\[
\prod_{e\subset \partial f} J_f
\subset
\prod_{e\subset \partial f} R.\]
We define values of $J$ on objects $p$ to be 
\[
\prod_{p\in \partial e;\, e\subset \partial f} 
J_f
\subset
\prod_{p\in \partial e;\,e\subset \partial f} R.\]

%We define a restriction map for $2$-morphisms: Let $f$ be a $2$-morphism in ${\cal{C}}$. 
%The $f:A\rightarrow \H\times\H$. Let $g$ be another $2$-morphism, 
%which we call a sub-$2$-morphism, which defines a sub-manifold of the image of $f$ so that
%$g=f|_B$, where $B\subset A$. 
%Then for each $2$-morphism $f$ and sub-$2$-morphisms $g$ , 
%we define a restriction morphism  \[J_f\rightarrow J_g\] 
%that sends each elements \[J_f(\Omega)\mapsto J_g(\Omega).\]

%Now, we define restriction map for $1$-morphisms: If $e$ is a $1$-morphism then it is a concatenation of finitely many geodesics. Let $e_1$ be a connected subset of the union of the geodesics of $e$. More precisely, let $e_1$ be a restriction of $e$ to a subinterval so that the boundary points of $e_1$ are cusps. Then there is a natural restriction morphism $J_e\rightarrow J_{e_1}$

%???
%sending $J_f(\Omega)$ from $J_e$ to $J_f(\Omega)$ in $J_{e_1}$. Note that this map is an inclusion, since every manifold with corners $Im(f)$ whose boundary contains the path $Im(e)$ also contains the sub-path $Im(e_1)$. 
%???

%!
%The statement in the question marks is true only if we assume that the values of $J_e$ are in $R$ but not in a countable product of $R$, $\prod R$. Then the restriction morphism is well-defined. It is actually an inclusion of sets.
%!

%Thus, $J$ is a contravariant functor from a $2$-category $C$ to a  $2$-category of sets.

The sheaf conditions for $1$-morphisms and the sheaf conditions for $2$-morphisms resemble the classical conditions for a presheaf to be a sheaf. 

Let $f_i:A_i\rightarrow [0,1]^2$ be a finite collection of disjoint $2$-morphisms, whose union is
 a morphism $f:A\rightarrow [0,1]^2$. 
We define a finite collection $f_{ij}^k$ of $1$-morphisms and $0$-morphisms (objects) such that the union 
$\bigcup_k Im(f^k_{ij})=Im(f_i)\cap Im(f_j)$ is a disjoint union of the intersection.

Then the equalizer 
\[J_f\rightarrow \prod_i J_{f_i}\rightrightarrows \prod_{ijk}J_{f_{ij}^k}\]
is exact (for a definition of equalizer one may consult \cite{Bor}).

Similarly, let $e$ be a $1$-morphism and let $\{e_i\}_i$ be a finite set of disjoint $1$-morphisms such that the union $\bigcup_i Im(e_i) =Im(e)$. We can write the intersection $Im(e_i)\cap Im(e_j)$ as a finite union of $0$-morphisms $\bigcup_k Im(e_{ij}^k)$, for some $1$-morphisms $e_{ij}^k$.

Then the equalizer 
\[J_e\rightarrow \prod_i J_{e_i}\rightrightarrows \prod_{ijk}J_{e_{ij}^k}\]
is exact.

The cochain is defined as
\[\prod_{
{\begin{small}
\begin{tabular}{c}
$p$: $0$-morph\\ 
\end{tabular}
\end{small}
}
}
J_p
\rightarrow 
\prod_{
{\begin{small}
\begin{tabular}{c}
$e$: $1$-morph\\ 
\end{tabular}
\end{small}
}
} J_e
\rightarrow
\prod_{
{\begin{small}
\begin{tabular}{c}
$f$: $2$-morph\\ 
\end{tabular}
\end{small}
}
}
J_f
\rightarrow
\prod_{
{\begin{small}
\begin{tabular}{c}
$g: \mbox{2-morph}$\\ 
$\partial g=\emptyset$
\end{tabular}
\end{small}
}
}J_g
\]
The maps $J_e\rightarrow J_p$ and $J_p\rightarrow J_f$ are surjective when they are defined, resembling flabby sheaves. Thus, we should have trivial $0$-th or $1$-st cohomology set. The only non-trivial cohomology will be $2$-nd cohomology set. The cocycle conditions for both non-commutative Hilbert modular symbols $c^1$ and $c^2$ can be interpreted as a particular case of maps 
\[\prod_{
{\begin{small}
\begin{tabular}{c}
$f$: $2$-morph\\ 
\end{tabular}
\end{small}
}
}
J_f
\rightarrow
\prod_{
{\begin{small}
\begin{tabular}{c}
$g: \mbox{2-morph}$\\ 
$\partial g=\emptyset$
\end{tabular}
\end{small}
}
}J_g
\]
For $c^1$, the boundary condition is that a union of two diangles with a common edge is a third diangle. One can think of the these three diangles as a boundary of a degenerate $3$-dimensional region. One can realize this cocycle condition as a sheaf-theoretic one by modifying the above definition so that the $2$-morphisms consists of a finite union of ideal diangles (without using the ideal triangles). Then the sheaf-theoretic $2$-nd cocycle condition is the one for non-commutative Hilbert modular symbol $c^1$.

If we are able to quotient  the $2$-category described in the beginning of this subsection by the $2$-morphisms generated by dianlges, then we have only two morphisms generated by ideal triangles.
The non-commutative Hilbert modular symbol $c^2$ is exactly the one that considers ideal triangles. Note that its cocycle relation for $c^2$ is satisfied up to $2$-morphisms generated by diangles.

\subsection{Explicit Computations. Multiple Dedekind Zeta Values}
In this Subsection, we make explicit computations of some ingredients in the non-commutative Hilbert modular symbol. In \cite{Man}, Manin compares explicit formulas of integrals in the non-commutative
modular symbol to multiple zeta values. The similarities are both in terms of infinite series formulas and in terms of formulas via iterated path integrals.
Here we compare certain integrals in the non-commutative Hilbert modular symbol to multiple Dedkeind zeta values (for multiple Dedekind zeta values, see \cite{MDZF}). Again the similarities are both in terms of infinite series formulas and 
in terms of formulas via iterated integrals over membranes.

We are going to consider the Fourier expansion of two Hilbert cusp forms $f$ and $g$. 
Let $\omega_f=fdz_1\wedge dz_2$, 
$\omega_g=gdz_1\wedge dz_2$ and 
$\omega_0=dz_1\wedge dz_2$. 
We are going to associate $L$-values to iterated integrals of the forms $\omega_f$ and $\omega_g$. 
The $L$-values would be iterated integrals over an union of diangles. One can think of a diangle connecting $0$ and $\infty$ as a segment or a real cone. The union will be a disjoint union of all such real cones connecting $0$ and $\infty$ or simply $\Im(\H)\times\Im(\H)$. We also recall the definition of a multiple Dedekind zeta values via (discrete) cone. Finally, we show analogous formulas for iterated $L$-values associated to Hilbert cusp forms and for multiple Dedekind zeta values.

We will be mostly interested in the modular symbol associated to a diangle. Let us recall what we mean by a diangle. 

Let $p_1,p_2,p_3,p_4$ be four cusp points. Let $\gamma_1\in GL_2(K)$ be a linear fractional transform that sends $\gamma_1(p_1)=0$,  $\gamma_1(p_2)=\infty$,  $\gamma_1(p_3)=1$. Let $\Delta$ be the image of the diagonal embedding of $\H^1$ into $\H^2$.
Then $0$, $1$ and $\infty$ are boundary points of $\Delta$. Let $\lambda(0,\infty)$ be the unique geodesic in $\Delta$ that connects $0$  and $\infty$. And let 
\[\lambda_1(p_1,p_2)=\gamma_1^{-1}\lambda(0,\infty)\]
be the pull-back of the geodesic $\lambda$ to a geodesic connecting $p_1$ and $p_2$.

Now consider the triple $p_1,p_2$ and $p_4$. Let $\gamma_2\in GL_2(K)$ be a linear fractional transform that sends $\gamma_2(p_1)=0$,  $\gamma_2(p_2)=\infty$ and $\gamma_2(p_4)=1$. 
Let $\Delta$ be the image of the diagonal embedding of $\H^1$ into $\H^2$.
Then $0$, $1$ and $\infty$ are boundary points of $\Delta$. Let $\lambda(0,\infty)$ be the unique geodesic in $\Delta$ that connects $0$  and $\infty$.
 And let 
\[\lambda_2(p_1,p_2)=\gamma_2^{-1}\lambda(0,\infty)\]
be the pull-back of the geodesic $\lambda$ to a geodesic connecting $p_1$ and $p_2$.

By a diangle, we mean a region in $\H^2\cup \P^1(K)$ of homotopy type of a disc, bounded by the geodesics $\lambda_1(0,\infty)$ and $\lambda_2(0,\infty)$. 

We are going to present a computation for the diangle $D_u$ defined by the points $(0,\infty,u^1,u^{-1}),$
where $u$ is a generator for the group of units modulo $\pm1$ in $K$. 
Let $(1)$ be the trivial permutation.
\begin{lemma} Let $u$ be a totally positive unit. Then
\[
\int\int_{D_u}^{(1)(1)}
e^
{2\pi i(\alpha_1z_1+\alpha_2z_2)}
dz_1\wedge dz_2
=
\frac{1}{(2\pi i)^2}\frac{u_2^2-u_1^2}{(\alpha_1u_1+\alpha_2u_2)(\alpha_1u_2+\alpha_2u_1)}\]
\end{lemma}
\proof Let $u_1$ and $u_2$ be the two embeddings of $u$ into $\R$.
Then $(0,\infty,u)$ can be mapped to $(0,\infty ,1)$ by
$\gamma_1=
\left(
\begin{tabular}{ll}
$u^{-1}$&$0$\\
$0$&$1$
\end{tabular}
\right).$
The geodesic $\lambda(0,\infty)$ can be parametrized by $(it,it)$ for $t\in \R$. Then the geodesic 
$\lambda_1(0,\infty)$ on the geodesic triangle $(0,\infty,u)$  can be parametrized by $\{(iu_1t,iu_2t)\,\,|\,\,t>0\}$. Similarly, the geodesic $\lambda_2(0,\infty)$ on the geodesic triangle $(0,\infty,u^{-1})$ can be parametrized by $\{(iu_2t,iu_1t)\,\,|\,\,t>0\}$. Then the diangle $D_u$ can be parametrized by
\[D_u=\{(z_1,z_2)\in \H^2\,|\,Re(z_1)=Re(z_2)=0, Im(z_1)\in\left(\frac{u_1}{u_2}t,\frac{u_2}{u_1}t\right), Im(z_2)=t\in (0,\infty)\}.\]
Then we have
\begin{align*}
\int\int_{D_u}^{(1)(1)}
e^
{2\pi i(\alpha_1z_1+\alpha_2z_2)}
dz_1\wedge dz_2
&=\int^0_\infty
\left(
\int_{\frac{u_2}{u_1}t}^{\frac{u_1}{u_2}t}
e^{2\pi i(\alpha_1z_1+\alpha_2 t)}
dz_1
\right)
dt=\\
&=\frac{1}{2\pi i\alpha_1}\int_\infty^0
\left(
e^{\alpha_1\frac{u_1}{u_2}t+\alpha_2t}
-e^{\alpha_1\frac{u_2}{u_1}t+\alpha_2t}
\right)
dt=\\
&=\frac{1}{(2\pi i)^2}\frac{1}{\alpha1}
\left(
\frac{1}{\alpha_1\frac{u_1}{u_2}+\alpha_2}
-\frac{1}{\alpha_1\frac{u_2}{u_1}+\alpha_2}
\right)=\\
&=\frac{1}{(2\pi i)^2}\frac{u_2^2-u_1^2}{(\alpha_1u_1+\alpha_2u_2)(\alpha_1u_2+\alpha_2u_1)}
\end{align*}
Therefore, one term of the Fourier expansion of a Hilbert cusp form paired with a symbol given by one diangle does not resemble a norm of an algebraic integer. However, if we integrate over an infinite union of diangles then a similarity with Dedekind zeta and with multiple Dedekind zeta occurs.

Consider the limit when $n\rightarrow \infty$ of $D_{u^n}$. It is the product of the two imaginary axes of the two upper half planes. Denote by 
\[
\Im(\H^2)=\Im(\H)\times\Im(\H).
\]
One can think of this region as an infinite union of diangles.

Denote by $\alpha z$ the sum of products $\alpha_1z_1+\alpha_2z_2$. 
Using the methods of \cite{MDZF}, Section 1, we obtain
\[
\frac{(2\pi i)^{-2}}{N(\alpha)N(\alpha+\beta)}
=\int_{\Im(\H^2)}^{(1)(1)} e^{2\pi i\alpha z}dz_1\wedge dz_2\cdot e^{2\pi i\beta z}dz_1\wedge dz_2\]
and
\[
\frac{1}{(2\pi i)^2}\frac{1}{N(\alpha)^3N(\alpha+\beta)^2}
=\int_{\Im(\H^2)}^{(1)(1)} e^{2\pi i\alpha z}dz_1\wedge dz_2\cdot (dz_1\wedge dz_2)\cdot (dz_1\wedge dz_2)
\cdot e^{2\pi i\beta z}dz_1\wedge dz_2\cdot (dz_1\wedge dz_2)\]

Let $f$ and $g$ be two cusp form of wights $(2k,2k)$ and $(2l,2l)$, respectively. 
Consider the Fourier expansion of both of the cusp forms.
Let 
\[
f=\sum_{\alpha>>0}a_\alpha e^{2\pi i \alpha z}
\]
and let 
\[
g=\sum_{\beta>>0}b_\beta e^{2\pi i \beta z}.
\]
Since $f$ is of weight $(2k,2k)$, we have that $a_{u\alpha}=a_{\alpha}$, where $u$ is a unit. For such a modular form the modular factor with respect to the transformation $z\rightarrow uz$ is $1$.
The $L$-values of $f$ is
\[
L_f(n)
=
\int_{\Im(\H^2)}^{(1)(1)}
\sum_{\alpha \in {\cal{O}}_K^+/U^+} 
a_\alpha e^{2\pi i\alpha z}
dz_1\wedge dz_2
\cdot 
(dz_1\wedge dz_2)^{\cdot (n-1)}
=
\frac{1}{(2\pi i)^{2n}}\sum_{\alpha\in {\cal{O}}_K^+/U^+}\frac{a_\alpha}{N(\alpha)^n}.
\]
Here ${\cal{O}}_K^+$ denotes the totally positive algebraic integers in $K$ and $U^+$ denotes the totally positive units. 

We recall some of the definitions from \cite{MDZF}. We fix a positive cone $C$ in ${\cal{O}}_K$, by which we mean
\[
C=\N\cup\{\alpha\in {\cal{O}}_K\,\,  |  \,\, a+b\epsilon,\, a,b\in\N\},
\]
where $\epsilon$ is a generator of the group of totally positive units.
By $\epsilon^k C$, we mean the collection of products $\epsilon^k \alpha$, where $\alpha$ varies in the cone $C$.

The following infinite sum is an example of a multiple Dedekind zeta value
\[
\zeta_{K;C,\epsilon^kC}(m,n)=\sum_{\alpha\in C}\sum_{\beta\in \epsilon^kC}
\frac{1}{N(\alpha)^mN(\alpha+\beta)^n}.
\]

Let $Z(m,n)=\sum_{k\in \Z}\zeta_{K;C,\epsilon^kC}(m,n),$ where $C$ is any set representing the totally positive algebraic integers ${\cal{O}}_K^+$ modulo totally positive units $U^+$.

\begin{lemma}
\label{lemma L}
The values $Z(m,n)$ are finite for $m>n>1$.
\end{lemma}
\proof Let $\epsilon$ be a generators of the group of totally positive units $U^+$ in $K$.
For the two real embeddings $\epsilon_1$ and $\epsilon_2$ of $\epsilon$, 
we can assume that $\epsilon_1>1>\epsilon_2$. Otherwise we can take its reciprocal. 

\begin{align}
\label{eq def}
Z(m,n)=&\sum_{k\in Z}\sum_{\alpha,\beta\in C}
\frac{1}{N(\alpha)^mN(\alpha+\epsilon^k \beta)^n}<\\
\label{eq e1}
<&\sum_{\alpha,\beta\in C} 
\frac{1}{N(\alpha)^m}\left(
\frac{1}{N(\alpha+\beta)^n}+\right.\\
\nonumber
&+\sum_{k=1}^\infty 
\frac{2^n}{\epsilon_1^k}
\left.\left(\frac{1}{\alpha_1^n\beta_2^n}+\frac{1}{\alpha_2^n\beta_1^n}\right)
\right)<\\
\nonumber
<&\sum_{\alpha,\beta\in C} 
\frac{1}{N(\alpha)^m}\left(
\frac{1}{N(\alpha+\beta)^n}+\right.\\
\label{eq a+b}
&+\sum_{k=1}^\infty 
\frac{2}{\epsilon_1^k}
\left.\left(\frac{N(\alpha+\beta)^n-N(\alpha)^n}{N(\alpha+\beta)^n}\right)
\right)=\\
\label{eq geometric}
=&\sum_{\alpha,\beta\in C} 
\frac{1}{N(\alpha)^m}\left(
\frac{1}{N(\alpha+\beta)^n}+\right.\\
\nonumber
&+\frac{2}{\epsilon_1-1}
\left.\left(\frac{N(\alpha+\beta)^n-N(\alpha)^n}{N(\alpha+\beta)^n}\right)
\right)=\\
\label{eq series}
=&\sum_{\alpha,\beta\in C} 
\frac{1}{N(\alpha)^mN(\alpha+\beta)^n}+\\
\nonumber
&+\frac{2}{\epsilon_1-1}\frac{1}{N(\alpha)^n}
-\frac{2}{N(\alpha)^{m-n}N(\alpha+\beta)^n}=\\
\label{eq result}
=&\zeta_K(C;m,n)-\frac{2}{\epsilon_1-1}\left(\zeta_K(C;n)+\zeta_K(C;m-n,n)\right).
\end{align}
Equation \eqref{eq def} is the definition. Inequality \eqref{eq e1} is based on the following:
 $\epsilon_2<1$ is replaced with $1$ then $k>0$. For $k<0$ we use $\epsilon_2^k=\epsilon_1^{-k}$.
We put $1$ for $\epsilon_1^{k}$ for $k<0$. The case $k=0$ is treated separately.
Finally we group the terms with equal powers of $\epsilon_1$. In the Inequality \eqref{eq a+b} we estimate the mixed terms in the brackets. In Equation \eqref{eq geometric} we take the sum of the geometric series in $\epsilon_1^{-1}$. Then in Equation \eqref{eq series} we open the brackets. And finally, in Equation \eqref{eq result}, we express the sums as a finite linear combinations of a Dedekind zeta value and multiple Dedekind zeta values. 

The following definition of an iterated $L$-value is the coefficient of one monomial from the non-commutative Hilbert modular symbol of type {\bf{b}}.
\begin{definition}
\label{def Lfg}
For a pair of  Hilbert cusp  forms $f$ and $g$ with Fourier expansion
\[
f=\sum_{\alpha>>0}a_\alpha e^{2\pi i \alpha z}
\text{ and }
g=\sum_{\beta>>0}b_\beta e^{2\pi i \beta z},
\]
we define the following iterated $L$-values
\begin{align*}
L_{f,g}(m,n)=&\int_{\Im(\H^2)}^{(1)(1)}\sum_{(\alpha,\beta)\in({\cal{O}}_K^+,{\cal{O}}_K^+)/U}
(a_\alpha e^{2\pi i\alpha z}dz_1\wedge dz_2)\cdot(dz_1\wedge dz_2)^{\cdot (m-1)}\cdot\\
&\cdot(b_\beta e^{2\pi i\beta z}dz_1\wedge dz_2)\cdot(dz_1\wedge dz_2)^{\cdot (n-1)}.
\end{align*}
\end{definition}

\begin{theorem}
\label{thm Lfg}
Using the above definition we have
\[L_{f,g}(m,n)=
\sum_{k\in\Z}\sum_{\alpha\in C, \beta\in \epsilon^kC}
\frac{a_\alpha b_\beta}{N(\alpha)^mN(\alpha+\beta)^n}\]
\end{theorem}
\proof
\begin{align*}
L_{f,g}(m,n)=&\int_{\Im(\H^2)}^{(1)(1)}\sum_{(\alpha,\beta)\in({\cal{O}}_K^+,{\cal{O}}_K^+)/U}
(a_\alpha e^{2\pi i\alpha z}dz_1\wedge dz_2)\cdot(dz_1\wedge dz_2)^{\cdot (m-1)}\cdot\\
&\cdot(b_\beta e^{2\pi i\beta z}dz_1\wedge dz_2)\cdot(dz_1\wedge dz_2)^{\cdot (n-1)}=\\
=&\sum_{(\alpha,\beta)\in({\cal{O}}_K^+,{\cal{O}}_K^+)/U}
\frac{a_\alpha b_\beta}{N(\alpha)^mN(\alpha+\beta)^n}=\\
=&\sum_{k\in\Z;\alpha,\beta\in C}
\frac{a_\alpha b_\beta}{N(\alpha)^mN(\alpha+\epsilon^k\beta)^n}=\\
=&\sum_{k\in\Z,}\sum_{\alpha\in C}\sum_{\beta\in \epsilon^kC}
\frac{a_\alpha b_\beta}{N(\alpha)^mN(\alpha+\beta)^n}.\\
\end{align*}
\qed

We would like to bring to the attention of the reader the Definition \ref{def Lfg} of the multiple $L$-
values. More specifically,  we would like to point out that the region of integration is an infinite union of 
diangles, (or equivalently an infinite union of real cones; see the beginning of this Section.) Note also 
that in  Theorem \ref{thm Lfg} the values of the multiple $L$-functoins are expressed as an infinite 
sums over different discrete cones, namely, over $\epsilon^kC$, for $k\in\Z$. However, a single real 
cone $D_u$ as in Lemma \ref{lemma L}, does not correspond to a single discrete cone. Only a good 
union of real cones $\Im(\H)\times\Im(\H)$ corresponds to a good union of discrete cones 
$\bigcup_{k\in\Z}(C,\epsilon^k C)$ as a fundamental domain of
$({\cal{O}}_K^+,{\cal{O}}_K^+)/U^+$.

\section*{Acknolwedgements} I would like to thank Ronnie Brown for his talk on higher categories at University of Durham, which helped me to define iterated integrals over membranes. I would like to thank Mladen Dimitrov for his interest in my approach to Hilbert modular symbols. Most of all I would like to thank to Yuri Manin for his inspiring talk at Max Planck Institute for Mathematics in Bonn, on non-commutative modular symbol and for his encouragement  and interest in the work presented here. Finally, I would like to thank the referees for the useful suggestions and commentaries. 

I would like to thank Max Planck Institute for Mathematics in Bonn and Marie Curie Research Training Network, Arithmetic Algebraic Geometry Network at Durham University for the financial support in 2004-2005 and in 2005-2006, respectively.

%%%%%%%%%%%%%%%%%%%%%%%%%%%%%%%%%%%%%%%%%%%

\renewcommand{\em}{\textrm}

\begin{small}

\renewcommand{\refname}{ {\flushleft\normalsize\bf{References}} }
    
\end{small}

\end{document}